\pdfoutput=1
\documentclass[twoside,12pt,a4paper]{cmft}

\pagespan{1}{?}

\newtheorem{theorem}{Theorem}[section]
\newtheorem{lemma}[theorem]{Lemma}

\newtheorem{alg}[theorem]{Algorithm}

\theoremstyle{remark}
\newtheorem{ex}[theorem]{Example}
\newtheorem{defn}[theorem]{Definition}
\newtheorem*{remark}{Remark}
\newtheorem*{acknowledgement}{Acknowledgment}

\usepackage{hyperref}
\usepackage{caption}
\usepackage{subfig}

\usepackage{tikz}
\usetikzlibrary{shapes,arrows}

\newcounter{minutes}
\divide\time by 60
\newcounter{hours}
\multiply\time by 60
\addtocounter{minutes}{-\time}

\newcommand{\C}{\mathbb{C}} 
\newcommand{\symD}{\Omega} 
\newcommand{\g}{\gamma} 

\newcommand{\symDisk}{\mathbb{D}}
\newcommand{\symQuad}{Q}
\newcommand{\symQuadC}{\tilde{Q}}
\newcommand{\symM}{\textrm{M}}

\begin{document}

\title[Conjugate Function Method for Numerical Conformal Mappings]{Conjugate Function Method for \\ Numerical Conformal Mappings}
\date{\today}

\author[H. Hakula]{Harri Hakula}
\email{harri.hakula@tkk.fi}
\address{Aalto University,
         Department of Mathematics and Systems Analysis,
         P.O. Box 11100,
         FI-00076 Aalto,
         Finland}
         
\author[T. Quach]{Tri Quach}
\email{tri.quach@tkk.fi}
\address{Aalto University,
         Department of Mathematics and Systems Analysis,
         P.O. Box 11100,
         FI-00076 Aalto,
         Finland}

\author[A. Rasila]{Antti Rasila}
\email{antti.rasila@iki.fi}
\address{Aalto University,
         Department of Mathematics and Systems Analysis,
         P.O. Box 11100,
         FI-00076 Aalto,
         Finland}



\begin{abstract}
\noindent We present a method for numerical computation of conformal mappings from simply or doubly connected domains onto so-called canonical domains, which in our case are rectangles or annuli. The method is based on conjugate harmonic functions and properties of quadrilaterals. Several numerical examples are given.
\end{abstract}

\maketitle



\section{Introduction}
In addition to their theoretical significance in complex analysis, 
conformal mappings are important in classical engineering applications, such as
electrostatics and aerodynamics \cite{SL}, but also in novel areas such as computer graphics and computational modeling \cite{bh,kss}. In this paper we study
numerical computation of conformal mappings $f$ of a domain $\symD
\subset \C$ into $\C$. We assume that the domain is bounded and that
there are either one or two simple (and non-intersecting) boundary
curves, i.e., the domain $\symD$ is either simply or doubly
connected. It is usually convenient to map the domains conformally
onto canonical domains, which are in our case rectangles $R_h = \{
z\in \C : 0< \mathrm{Re}\,z <1, \, 0< \mathrm{Im}\,z <h \}$ or
annuli $A_r = \{z \in \C: e^{-r}<|z|<1\}$. While the existence of
such conformal mappings is expected because of Riemann's mapping
theorem, it is usually  not possible to obtain a formula or other
representation for the mapping analytically.


Several different algorithms for numerical computation of conformal
mappings have been described in the literature. One popular method
involves the Schwarz-Christoffel formula, which can also be generalized
for doubly connected domains. A widely used MATLAB
implementation of this method is due to Driscoll \cite{Dri} and
a FORTRAN version for the doubly connected case is due to Hu \cite{Hu}. For
theoretical background concerning these methods see
\cite{DT,DrVa,Tre}. In addition, there are several
approaches which do not involve the Schwarz-Christoffel formula, e.g.,
the Zipper algorithm of Marshall \cite{Mar, MR}. A method involving the harmonic conjugate function is presented in \cite[pp. 371-374]{Hen3}, but this method is different from ours as it does not use quadrilaterals.
For an overview of numerical conformal mappings and moduli of
quadrilaterals, see \cite{PS}. Historical remarks and an outline of
development of numerical methods in conformal mappings is given in
\cite{DT,Kyt,Porter}. 

In this paper, we
present a new method for constructing numerical conformal mappings. The method is based on the harmonic conjugate
function and properties of quadrilaterals, which together form the 
foundation of our numerical algorithm. The algorithm is based on
solving numerically the Laplace equation subject to Dirichlet-Neumann 
mixed-type boundary conditions, which is described in \cite{HRV}. To the best of our knowledge, this is the first attempt to construct conformal mappings by using $hp$-FEM. It should be noted, that the presented method is not restricted to polygonal domains, and can be used with domains with curvilinear boundary as well. 

The outline of the paper is as follows.
First the preliminary concepts are introduced and then the new algorithm
is described in detail. Before the numerical examples, the
computational complexity and some details of our implementation are
discussed.
The numerical examples are divided into three sections:
validation against the Schwarz-Christoffel toolbox, the analytic example, simply connected domains, and finally ring domains.

\section{Foundations of the Conjugate Function Method} \label{sec: conj-func}

In this section we introduce the required concepts from function
theory, and present a proof of a fundamental result leading to a numerical algorithm.

\begin{defn} (Modulus of a Quadrilateral) \\
A Jordan domain $\symD$ in $\C$ with marked (positively ordered) points 
$z_1,z_2,z_3,z_4\in \partial \symD$ is called a {\it quadrilateral}, and is denoted 
by $\symQuad = (\symD;z_1,z_2,z_3,z_4)$. Then there is a canonical conformal map of the quadrilateral $\symQuad$ onto a rectangle $R_h = (\symD';1+ih,ih,0,1)$, with the vertices 
corresponding, where the quantity $h$ defines the  {\it modulus of a quadrilateral}
$\symQuad$. We write
\[
\symM(\symQuad) = h.
\]
\end{defn}
Note that the modulus $h$ is unique.
\begin{defn} (Reciprocal Identity) \\
 It is clear by the geometry
\cite[p. 15]{LV} or \cite[pp. 53-54]{PS} that the following reciprocal
identity holds:
\begin{equation} \label{eqn: recip}
\symM(\symQuad)\, \symM(\symQuadC) =1,
\end{equation}
where $\symQuadC= (\symD; z_2,z_3,z_4, z_1)$ is called the
 {\it conjugate quadrilateral} of $\symQuad$.
\end{defn}
 For basic properties of modulus of quadrilaterals, we refer the reader to \cite{LV} and \cite[Chapter~2]{PS}.

\begin{remark}
The identity \eqref{eqn: recip} leads to a method for estimating the numerical accuracy of the modulus. For discussion and several numerical examples, see \cite{HRV}.
\end{remark}




\subsection{Dirichlet-Neumann Problem}
It is well known that one can express the modulus of a quadrilateral $Q$ in terms of the solution of the Dirichlet-Neumann mixed boundary value problem \cite[p. 431]{Hen3}.

Let $\symD$ be a domain in the complex plane whose boundary $\partial \symD$ consists of a finite number of regular Jordan curves, so that at every point, except possibly at finitely many points of the boundary, a normal is defined.  Let $\partial \symD =A \cup B$ where $A, B$ both are unions of regular Jordan arcs such that $A \cap B$ is finite. Let $\psi_A$, $\psi_B$ be real-valued continuous functions defined on $A, B$, respectively. Find a function $u$ satisfying the following conditions:
\begin{enumerate}
\item $u$ is continuous and differentiable in $\overline{\symD}$.
\item $u(t) = \psi_A(t),\qquad \textrm{for all } \, t \in A$.
\item If $\partial/\partial n$ denotes differentiation in
the direction of the exterior normal, then
\[
\frac{\partial}{\partial n} u(t)=\psi_B(t),\qquad \textrm{for all } \, t \in  B.
\]
\end{enumerate}
The problem associated with the conjugate quadrilateral $\symQuadC$ is called the {\it conjugate Dirichlet-Neumann problem}.

Let $\gamma_j, j=1,2,3,4$ be the arcs of $\partial \symD$ between $(z_1, z_2)\,,$ $(z_2, z_3)\,,$ $(z_3, z_4)\,,$ $(z_4, z_1),$ respectively. Suppose that $u$ is the (unique) harmonic solution of the Dirichlet-Neumann problem with mixed boundary values of $u$ equal to $0$ on $\gamma_2$, equal to $1$ on $\gamma_4$, and $\partial u/\partial n = 0$ on $\gamma_1, \gamma_3$. Then by \cite[Theorem 4.5]{Ahl} or \cite[Theorem 2.3.3]{PS}:
\begin{equation} \label{qmod}
\symM(\symQuad)= \iint_\symD |\nabla  u|^2\,dx \, dy.
\end{equation}

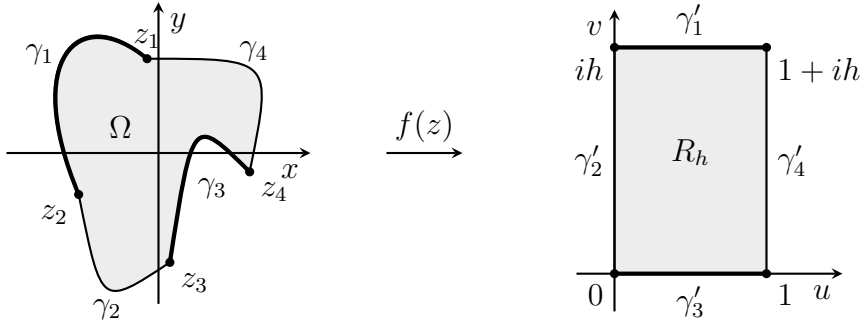
\begin{figure}[!ht]
\begin{center}
\begin{tikzpicture}[scale=1,>=stealth]
\draw[rounded corners=2pt,thick,fill=black!7!white,xshift=-0.3cm] (0.45,-0.45) coordinate (z_3) .. controls (0.75,1.5) .. (1.5,0.75) coordinate (z_4) {node[pos=0.3,anchor=north west]{$\gamma_3$}} .. controls (1.8,2.25) .. (0.15,2.25) coordinate (z_1) {node[pos=0.5,anchor=south]{$\gamma_4$}} .. controls (-0.75,3) and (-1.5,2.25) .. (-0.75,0.45) coordinate (z_2) {node[pos=0.5,anchor=east]{$\gamma_1$}} .. controls (-0.45,-1.05) .. (0.45,-0.45) {node[pos=0.5,anchor=north]{$\gamma_2$}};
\draw[rounded corners=2pt,ultra thick,xshift=-0.3cm] (z_1) .. controls (-0.75,3) and (-1.5,2.25) .. (z_2);
\draw[rounded corners=2pt,ultra thick,xshift=-0.3cm] (z_3) .. controls (0.75,1.5) .. (z_4);
\filldraw (z_1) circle (1.5pt) node[anchor=south]{$z_1$};
\filldraw (z_2) circle (1.5pt) node[anchor=north east]{$z_2$};
\filldraw (z_3) circle (1.5pt) node[anchor=north west]{$z_3$};
\filldraw (z_4) circle (1.5pt) node[anchor=north west]{$z_4$};
\draw (-0.5,1.2) node[anchor=base]{$\symD$};

	\draw[->,thick,yshift=1cm] (0,-2) -- (0,2) node [anchor=north west] {$y$};
	\draw[->,thick,yshift=1cm] (-2,0) -- (2,0) node [anchor=north east] {$x$};

\draw[thick,fill=black!7!white,xshift=6cm,yshift=-0.6cm] (0,0) rectangle (2,3);
\draw[ultra thick,xshift=6cm,yshift=-0.6cm] (0,0) coordinate (w_1) -- (2,0) coordinate (w_2) {node[pos=0.5,anchor=north]{$\gamma_3'$}} (2,3) coordinate (w_3) -- (0,3)  coordinate (w_4) {node[pos=0.5,anchor=south]{$\gamma_1'$}};
\draw[thin,xshift=6cm,yshift=-0.6cm] (w_2) -- (w_3) {node[pos=0.5,anchor=west]{$\gamma_4'$}};
\draw[thin,xshift=6cm,yshift=-0.6cm] (w_4) -- (w_1) {node[pos=0.5,anchor=east]{$\gamma_2'$}};

\filldraw (w_1) circle (1.5pt) node[anchor=north east]{$0$};
\filldraw (w_2) circle (1.5pt) node[anchor=north west]{$1$};
\filldraw (w_3) circle (1.5pt) node[anchor=north west]{$1+ih$};
\filldraw (w_4) circle (1.5pt) node[anchor=north east]{$ih$};

	\draw[->,thick,xshift=6cm,yshift=-0.6cm] (0,-.5) -- (0,3.5) node [anchor=north east] {$v$};
	\draw[->,thick,xshift=6cm,yshift=-0.6cm] (-.5,0) -- (3,0) node [anchor=north east] {$u$};
\draw[xshift=6cm,yshift=-0.5cm] (1,1.4) node[anchor=base]{$R_h$};
\draw[->,thick,xshift=3cm,yshift=1cm] (0,0) -- (1,0) node [midway,anchor=south]{$f(z)$};
\end{tikzpicture}
\caption{Dirichlet-Neumann boundary value problem. Dirichlet and Neumann boundary conditions are marked with thin and thick lines, respectively. \label{fig: module_DN}}
\end{center}
\end{figure}

Suppose that $\symQuad$ is a quadrilateral, and $u$ is the harmonic solution of the Dirichlet-Neumann problem and let $v$ be a conjugate harmonic function of $u$ such that $v(\textrm{Re}\, z_3, \textrm{Im}\, z_3) = 0$. Then $f = u + iv$ is an analytic function, and it maps $\symD$ onto a rectangle $R_h$ such that the image of the points $z_1,z_2,z_3,z_4$ are $1+ih, ih,0,1$, respectively. Furthermore by Carath\'{e}odory's theorem \cite[Theorem 5.1.1]{Kra}, $f$ has a continuous boundary extension which maps the boundary curves $\g_1, \g_2, \g_3, \g_4$ onto the line segments $\g_1', \g_2', \g_3', \g_4'$, see Figure \ref{fig: module_DN}.

\begin{lemma} \label{lemma: conj-h2}
Let $\symQuad$ be a quadrilateral with modulus $h$, and let $u$ be the harmonic solution of the Dirichlet-Neumann problem. Suppose that $v$ is the harmonic conjugate function of $u$, with $v({\rm Re}\, z_3, {\rm Im}\, z_3) = 0$. If  $\tilde{u}$ is the harmonic solution of the Dirichlet-Neumann problem associated with the conjugate quadrilateral $\symQuadC$, then $v = h\tilde{u}$.
\end{lemma}

\proof It is clear that $v, \tilde{u}$ are harmonic. Thus $\tilde{v} = h\tilde{u}$ is harmonic, and $v$ and $\tilde{v}$ are both constant on $\gamma_1, \gamma_3$. By Cauchy-Riemann equations, we obtain $\langle \nabla u, \nabla v \rangle = 0$. We may assume that the gradient of $u$ does not vanish on $\gamma_2, \gamma_4$. In particular, on $\gamma_4$, we have $n = \nabla u / |\nabla u|$, where $n$ is the exterior normal of the boundary. On the other hand, on $\gamma_2$, we have $n = -\nabla u / |\nabla u|$. Therefore, we have 
\[
\frac{\partial v}{\partial n} = \langle \nabla v, n \rangle = \pm\frac{1}{|\nabla u|} \, \langle \nabla v, \nabla u \rangle = 0.
\]
By the definition of $\tilde{u}$, we get
\[
\frac{\partial \tilde{v}}{\partial n}  = h\frac{\partial \tilde{u}}{\partial n} = 0,
\]
on $\gamma_2, \gamma_4$. Thus $v$ and $\tilde{v}$ satisfy the same boundary conditions on $\gamma_2, \gamma_4$. Then by \eqref{eqn: recip} and the uniqueness theorem for harmonic functions \cite[p. 166]{Ahl2}, we conclude that $v = \tilde{v}$. \qed

Suppose that $f=u+iv$, where $u$ and $v$ are as in Lemma \ref{lemma: conj-h2}. Then it is easy to see that the image of equipotential curves of the functions $u$ and $v$ are parallel to the imaginary and the real axis, respectively.

Finally, we note that the function $f$ constructed this way is univalent. To see this, suppose that $f$ is not univalent. Then there exists points $z_1, z_2 \in \symD$ and $z_1 \not= z_2$ such that $f(z_1) = f(z_2)$. Thus $\mathrm{Re}\, f(z_1) = \mathrm{Re}\, f(z_2)$, so $z_1$ and $z_2$ are on the same equipotential curve $C$ of $u$. Similarly for imaginary part, $z_1$ and $z_2$ are on the same equipotential curve $\tilde{C}$ of $v$. Then by the above fact of equipotential curves, it follows that $z_1 = z_2$, which is a contradiction.

\subsection{Ring Domains}
Let $E$ and $F$ be two disjoint and connected compact sets in the extended 
complex plane ${\C_\infty} = \C \cup \{\infty\}$. Then one of the sets $E, F$ is bounded and
without loss of generality we may assume that it is $E$. Then a set $R={\C_\infty} \backslash (E \cup F)$ is connected and is called a {\it ring domain}. The {\it capacity} of $R$ is defined by
\[
\textrm{cap} R = \inf_u \iint_R |\nabla u|^2 \, dx \, dy,
\]
where the infimum is taken over all non-negative, piecewise differentiable functions $u$ with compact support in $R\cup E$ such that $u=1$ on $E$. Suppose that a function $u$ is defined on $R$ with $1$ on $E$ and $0$ on $F$. Then if $u$ is harmonic, it is unique and it minimizes the above integral. The conformal modulus of a ring domain $R$ is defined by $\symM(R) = 2\pi / \textrm{cap} R$. The ring domain $R$ can be mapped conformally onto the annulus $A_r$, where $r = \symM(R)$. In \cite{BSV} numerical computation of modulus of several ring domains is studied.

\section{Conjugate Function Method}

Our aim is to construct a conformal mapping from a quadrilateral
$\symQuad = (\symD; z_1, z_2, z_3, z_4)$ onto a rectangle $R_h$, where
$h$ is the modulus of the quadrilateral $\symQuad$. Here the points
$z_j$ will be mapped onto the corners of the rectangle $R_h$. 
In the standard algorithm the required steps are the following:
\begin{alg}(Conformal Mapping) \label{alg: conf}
\\
\vspace*{-0.7cm}
\begin{enumerate}
\item  Find a harmonic solution for a Dirichlet-Neumann problem associated with a quadrilateral.
\item Solve the Cauchy-Riemann differential equations in order to obtain an analytic function that maps our domain of interest onto a rectangle.
\end{enumerate}
\end{alg}

The Dirichlet-Neumann problem can be solved by using any suitable numerical method.
One could also solve the Cauchy-Riemann equations numerically (see e.g. \cite{Bra}) but it is not necessary.
Instead we solve $v$ directly from 
the conjugate problem, which is usually computationally much more
efficient, because the mesh and the discretized system used in solving
the potential function $u$ 
can be used for solving $v$ as well.

This new algorithm is as follows:
\begin{alg}(Conjugate Function Method) \label{alg: hqr}
\\
\vspace*{-0.7cm}
\begin{enumerate}
\item Solve the Dirichlet-Neumann problem to obtain $u_1$ and compute the modulus $h$.
\item Solve the Dirichlet-Neumann problem associated with $\symQuadC$ to obtain $u_2$.
\item Then $f = u_1 + ihu_2$ is the conformal mapping from $\symQuad$ onto $R_h$ such that the vertices $(z_1,z_2,z_3,z_4)$ are mapped onto the corners $(1+ih,ih,0,1)$.
\end{enumerate}
\end{alg}

\begin{figure}
\centering
\subfloat[Ring domain with Dirichlet data 0, and 1, on the outer and
inner boundaries, respectively.]{\parbox{.45\textwidth}{\centering\includegraphics[width=.35\textwidth]{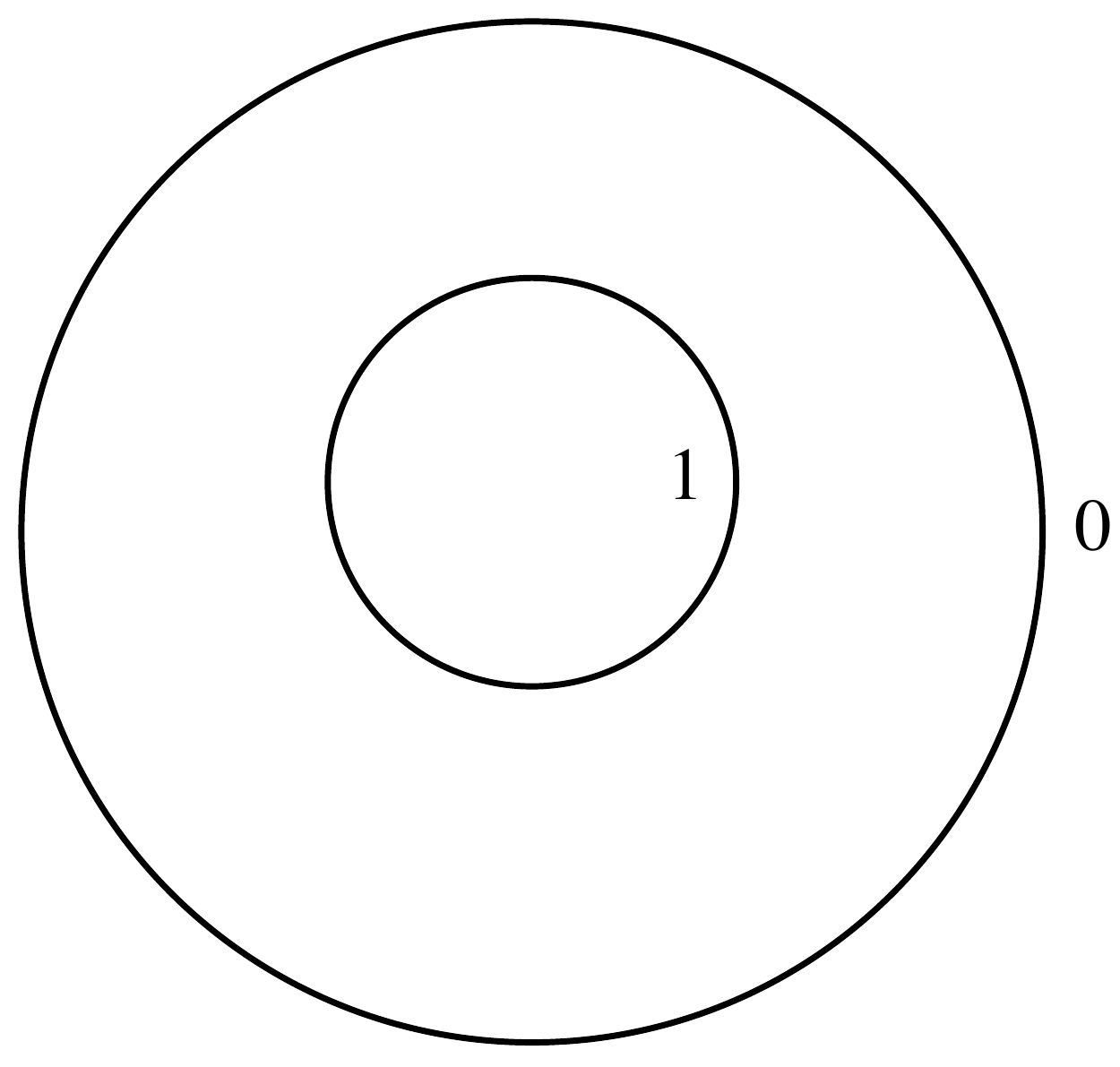}}}\hfill
\subfloat[Ring domain: Solution of the Dirichlet problem with contour lines.]{\parbox{.45\textwidth}{\centering\includegraphics[width=.35\textwidth]{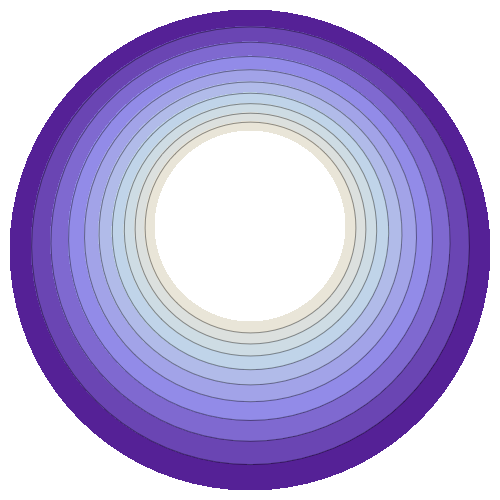}}}\\
\subfloat[Conjugate problem for the cut domain with new Dirichlet data along the both sides of
the cut.]{\parbox{.45\textwidth}{\centering\includegraphics[width=.35\textwidth]{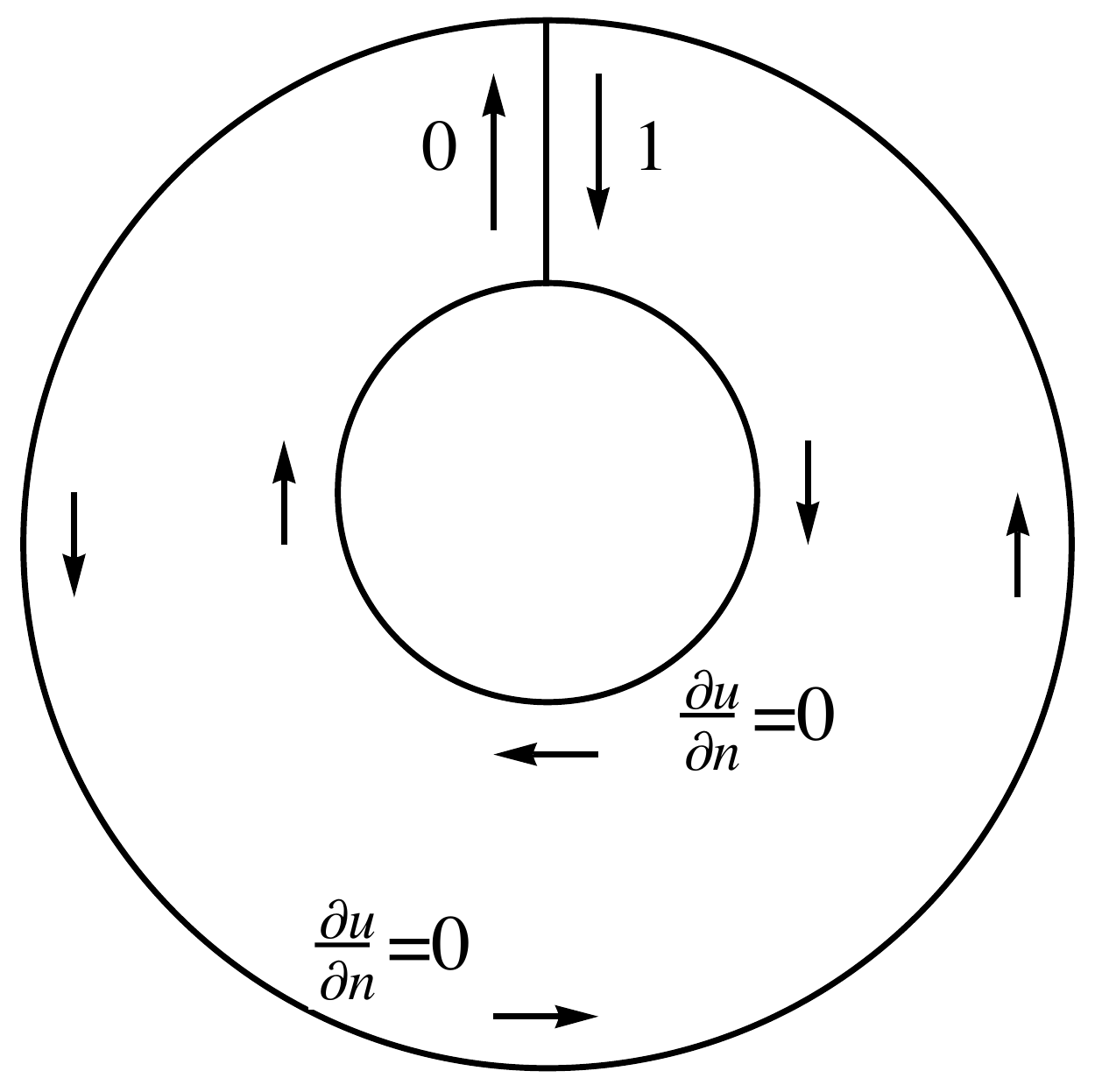}}}\hfill
\subfloat[Conjugate problem: Solution of the conjugate problem  with contour lines.]{\parbox{.45\textwidth}{\centering\includegraphics[width=.35\textwidth]{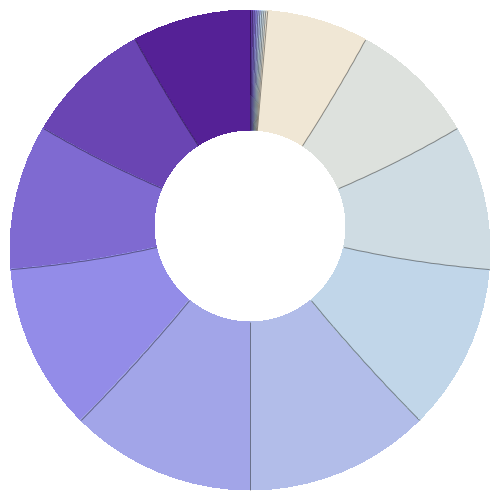}}}\\
\subfloat[Mapped annulus.]{\parbox{.5\textwidth}{\centering\includegraphics[width=.35\textwidth]{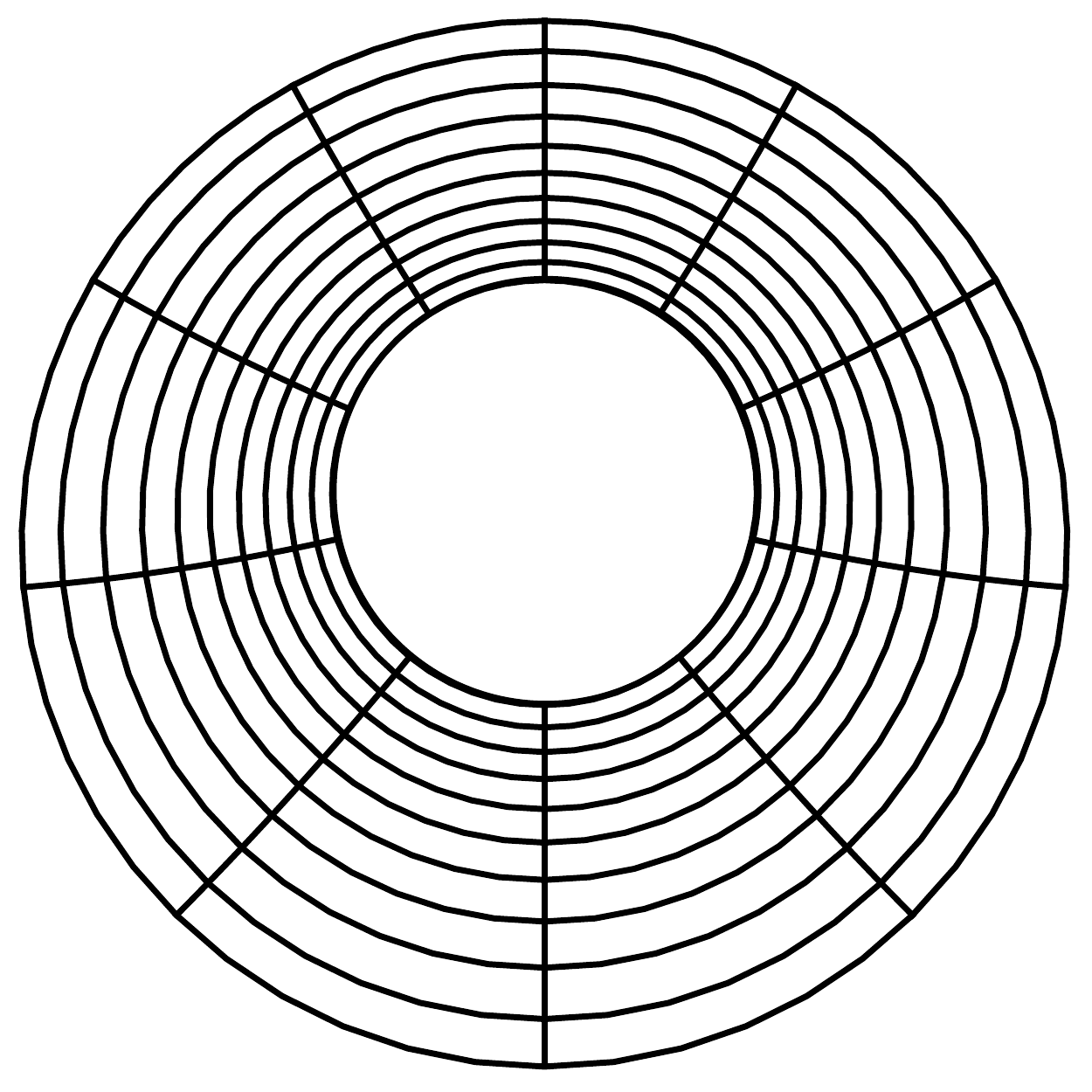}}}
\caption{Conjugate Function Method for Ring Domains.}\label{fig:algorithm}
\end{figure}

 In the case of ring domains, the construction of the conformal mapping is
 slightly different. The necessary steps are described below and in
 Figure \ref{fig:algorithm}.

\begin{alg}(Conjugate Function Method for Ring Domains) \label{alg: hqrring}
\\
\vspace*{-0.7cm}
\begin{enumerate}
\item Solve the  Dirichlet problem to obtain the potential function
$u$ and the modulus $\symM(R)$.
\item  Cut the ring domain through the steepest descent curve which
is given by the gradient of the potential function $u$ and obtain a quadrilateral where the Neumann condition is on the
steepest descent curve and the Dirichlet boundaries remain as before. 
\item Use the Algorithm \ref{alg: hqr}.
\end{enumerate}
\end{alg}
Note that the choice of the steepest descent curve is not unique due
to the implicit orthogonality condition.

\section{Implementation Aspects}
\label{sec:implementation}
The $hp$-FEM implementation we are using is described in detail in \cite{HRV}.
For elliptic problems, the superior accuracy of the higher order methods with
relatively small number of unknowns
has to be balanced against the much higher integration cost
and the cost of evaluating the solution at any given point in the
domain. It should be emphasized though, that the conjugate function method
is not dependent on any particular numerical PDE solution technique.
Indeed, we fully expect that similar results could be obtained
with, for instance, fine-tuned integral equation solvers. 

In the context of solution of the conjugate pair problems, it is
obvious that we have to integrate only once. Moreover, the
factorization of the resulting discretized systems can be, for the most
part, used in both problems without any extra work. Therefore, although in
principle two problems are solved, in practice the work is almost proportional
to that of one.

However, the computation of the contour lines necessarily involves a
large number of evaluations of the solution, that also become more
expensive as the order of the method increases. 

\subsection{$\boldsymbol{hp}$ -FEM}

Here we give a short overview of the $hp$-FEM following closely to the one in \cite{HRV}.
In the $p$-version of the FEM the polynomial degree $p$ is used to control
the accuracy of the solution while keeping the mesh fixed in contrast to the
$h$-version or standard finite element method, where the polynomial degree is
constant but the mesh size varies.
Thus, the $p$-version is often referred to as the $p$-extension. 
The $hp$-method simply combines the $h$- and $p$- refinements.

These different refinement strategies also imply different sets of
unknowns or degrees of freedom:
In the $h$-version or the standard finite element method, the unknowns or degrees
of freedom are associated with values at specified locations of the discretization
of the computational domain, that is, the nodes of the mesh.
In the $p$-method, the unknowns are coefficients of some polynomials that
are associated with topological entities of the elements, nodes, sides, and the interior.

For optimal $hp$-convergence one should refine the mesh geometrically toward
corners and let the degree of polynomial shape functions increase with distance from the corners.
For an example of such a mesh see Figure~\ref{fig: mesh}. In the examples
computed below, we have used a constant value of the order over the whole mesh.

In the following one way to construct a $p$-type quadrilateral element is given.
The construction of triangles follows similar lines. First of all, the choice of
shape functions is not unique.
We use the so-called hierarchic integrated Legendre shape functions.

Legendre polynomials of degree $n$ can be defined by a recursion formula
\[
(n+1)P_{n+1}(x) - (2n+1) x P_n(x) + n P_{n-1}(x) =0,
\]
where $P_0(x)=1$ and $P_1(x) = x$.

The derivatives can similarly be computed by using the recursion
\[
(1-x^2)P'_n(x)=-n x P_n(x) + n P_{n-1}(x).
\]

The integrated Legendre
polynomials are defined for $x\in[-1,1]$ as
\[
\phi_n(\xi)=\sqrt{\frac{2n-1}{2}}\int_{-1}^{\xi}P_{n-1}(t)\,dt,\quad n=2,3,\ldots ,
\]
and can be rewritten as linear combinations of Legendre polynomials
\[
\phi_n(\xi)=
\frac{1}{\sqrt{2(2n -1)}}\left(P_n(\xi) - P_{n-2}(\xi)\right),\quad n=2,3,\ldots .
\]
The normalizing coefficients are chosen so that
\[
\int_{-1}^{1}\frac{d \phi_i(\xi)}{d\xi} \frac{d \phi_j(\xi)}{d\xi} \,d\xi = 
	\delta_{ij}, \quad i,j \geq 2.
\]

Using these polynomials
we can now define the shape functions for a quadrilateral reference element
over the domain $[-1,1]\times[-1,1]$. The shape functions are divided into three categories:
nodal shape functions, side modes, and internal modes.

There are four nodal shape functions.
\begin{align*}
N_1(\xi,\eta) &= \frac{1}{4}(1-\xi)(1-\eta), \quad
N_2(\xi,\eta) = \frac{1}{4}(1+\xi)(1-\eta), \\
N_3(\xi,\eta) &= \frac{1}{4}(1+\xi)(1+\eta), \quad
N_4(\xi,\eta) = \frac{1}{4}(1-\xi)(1+\eta), 
\end{align*}
which taken alone define the standard four-node quadrilateral 
finite element. There are $4(p-1)$ side modes associated with
the sides of a quadrilateral $(p\geq 2)$, with $i=2,\ldots,p$,
\begin{align*}
N_i^{(1)}(\xi,\eta) &= \frac{1}{2} (1-\eta) \phi_i(\xi), \quad
N_i^{(2)}(\xi,\eta) = \frac{1}{2} (1+\xi) \phi_i(\eta),\\
N_i^{(3)}(\xi,\eta) &= \frac{1}{2} (1+\eta) \phi_i(\eta),\quad
N_i^{(4)}(\xi,\eta) = \frac{1}{2} (1-\xi) \phi_i(\xi).
\end{align*}
For the internal modes we choose the $(p-1)(p-1)$ shape functions
\[
N_{i,j}^0(\xi,\eta)=\phi_i(\xi)\phi_j(\eta),\quad i=2,\ldots,p,\quad j=2,\ldots,p.
\]
The internal shape functions are often referred to as bubble-functions.

Note that some additional book-keeping is necessary.
The Legendre polynomials have the property $P_n(-x)=(-1)^n P_n(x)$.
This means that every edge must be globally parameterized in the same way
in both elements to where it belongs.
Otherwise unexpected cancellation in the degrees of freedom associated
with the odd edge modes could occur.

\subsection{Solution of Linear Systems}
Let us divide the degrees of freedom of a discretized quadrilateral into five sets, internal
and boundary degrees of freedom. The sets are denoted
$B, D_0, D_1, N^0$, and $N^1$, for internal, Dirichlet $u=0$,
Dirichlet $u=1$, 
Neumann with Dirichlet $u=0$ in the conjugate problem, and 
Neumann with Dirichlet $u=1$ in the conjugate problem,  degrees of
freedom, respectively.

The discretized system is
\[
A =
\begin{pmatrix}
A_{BB}     & A_{B N^1} & A_{B N^0} & A_{B D_1} & A_{B D_0} \\
A_{N^1 B} & A_{N^1 N^1} & A_{N^1 N^0} & A_{N^1 D_1} & A_{N^1 D_0} \\
A_{N^0 B} & A_{N^0 N^1} & A_{N^0 N^0} & A_{N^0 D_1} & A_{N^0 D_0} \\
A_{D_1 B} & A_{D_1 N^1} & A_{D_1 N^0} & A_{D_1 D_1} & A_{D_1 D_0} \\
A_{D_0 B} & A_{D_0 N^1}  & A_{D_0 N^0} & A_{D_0 D_1} & A_{D_0 D_0} \\
\end{pmatrix}.
\]
Taking the Dirichlet boundary conditions into account, we arrive at
the following system of equations, using $x_{D_1} = \mathbf{1}$,
\[
\begin{pmatrix}
A_{BB}     & A_{B N^1} & A_{B N^0} \\
A_{N^1 B} & A_{N^1 N^1} & A_{N^1 N^0} \\
A_{N^0 B} & A_{N^0 N^1} & A_{N^0 N^0} \\
\end{pmatrix} 
\begin{pmatrix}
x_{B}  \\
x_{N^1} \\
x_{N^0} \\
\end{pmatrix} =
-\begin{pmatrix}
A_{B D_1}  \mathbf{1}   \\
A_{N^1 D_1}   \mathbf{1} \\
A_{N^0 D_1}  \mathbf{1} \\
\end{pmatrix}.
\]
For the conjugate problem, simply change the roles of $D_1\leftrightarrow N^1$
and $D_0\leftrightarrow N^0$. Note that $A_{BB} $ is present in both systems
and thus \textit{has to be factored only once}.

\subsection{Evaluation of Contour Lines}

Let $u$ and $v$ be solutions of the Dirichlet-Neumann problem and its conjugate problem, respectively.
In computing the contour lines, the solutions and their gradients have
to be evaluated at many points $(x,y)$.
Evaluation of the solution in $hp$-FEM is more complicated than in the
standard FEM. In a reference-element-based system such as ours, in
order to evaluate the solution at point $(x,y)$ one must first find
the enclosing element and then the local coordinates of the point on
that element. Then every shape function has to be  evaluated at the local
coordinates of the point. This is outlined in detail in Algorithm \ref{alg:hpeval}.
Similarly evaluation of the gradient of the solution requires two
polynomial evaluations per one geometric search.

\begin{alg}(Evaluation of $u(x,y)$) \label{alg:hpeval}
\\
\vspace*{-0.7cm}
\begin{enumerate}
\item Find the enclosing element of $(x,y)$.
\item Find the local coordinates $(\xi,\eta)$ on the element.
\item Evaluate the shape functions $\phi_i(\xi,\eta)$.
\item Compute the linear combination of the shape functions $\sum_i
  c_i \phi_i(\xi,\eta)$, where $c_i$ are the coefficients from the
  solution vector.
\end{enumerate}
\end{alg}

Finding the images of the canonical domains is equivalent to finding
the corresponding contour lines of $u$ and $v$. 
Since both solutions have been computed on the same mesh,
evaluating the solutions and their gradients at the same point is
straightforward. In Algorithm \ref{alg:contour} the two-level line
search is described in detail.

\begin{alg}(Tracing of Contour Lines: $u(x,y) = c = \text{const}$.) \label{alg:contour}
\\
\vspace*{-0.7cm}
\begin{enumerate}
\item Find the solutions $u(x,y)$ and $v(x,y)$.
\item Set the step size $\sigma$ and the tolerance $\epsilon$.
\item Choose the potential $c$.
\item Search along the Neumann boundary for the point $(x,y)$ such
  that $u(x,y)~=~c$.
\item Take a step of length $\sigma$ along the contour line of $u(x,y)$ in the direction of
  $\nabla v(x,y)$ to a new point $(\hat{x},\hat{y})$.
\item Correct the point $(\hat{x},\hat{y})$  by searching in the orthogonal direction,
  i.e., $\nabla u(\hat{x},\hat{y})$, until $|u(\hat{x},\hat{y}) - c| < \epsilon$  is achieved.
\item Set $(x,y) = (\hat{x},\hat{y})$ and repeat until the opposite Neumann boundary has been reached.
\end{enumerate}
\end{alg}


\subsection{On Computational Complexity}

The solution time of a single problem can be divided into three parts,
the setting up of the problem, the solution of the Dirichlet-Neumann and
its conjugate problem, and the evaluation of the mappings.
In short, in the absence of
fully automatic mesh generators for this class of problems, the setting up of
the problem remains the most time consuming part.
We have implemented the algorithm using Mathematica~8.

The time and memory requirements, in terms of degrees of freedom,
have been reported for the Dirichlet-Neumann problems in \cite{HRV}.  
In the examples below, the solution times vary from few seconds to at most
two minutes on standard hardware (as defined in Table \ref{table: countour}).
It should be noted that for comparable accuracy on polygonal domains, the Schwarz-Christoffel toolbox
is superior to our implementation. 

\begin{table}[h]
\caption{Effect of $p$ and $\epsilon$ on contour lines computations. 
In geometry of Figure \ref{fig: disk-in-pentagon}, ten contours of
(radial) $u(x,y) = c_u$ and (circular) $v(x,y) = c_v$ have been computed with $\sigma = 1/4$.
Times are normalized so that for $p$=8, $\epsilon$ = 1/100,  time $=1$. 
The time units are thus 1s and 46s for radial and circular contours, respectively.
  (Apple Mac Pro 2009 Edition 2.26 GHz, Mathematica 8.0.4)
  }\label{table: countour}
\begin{tabular}{|ll|cc|}
\hline
$p$ & $\epsilon$ & time for $u(x,y) = c_u$ & time for $v(x,y) = c_v$ \\ \hline
4 & 1/100    & 0.44 & 0.43 \\
4 & 1/1000   & 0.41 & 0.77 \\
4 & 1/10000  & 1.51 & 1.19 \\
8 & 1/100    & 1    & 1    \\
8 & 1/1000   & 1.00 & 1.82 \\
8 & 1/10000  & 0.99 & 2.66 \\
12 & 1/100   & 2.29 & 2.31 \\
12 & 1/1000  & 2.26 & 4.16 \\
12 & 1/10000 & 2.26 & 6.07 \\
\hline
\end{tabular}
\end{table}
Estimating the computational complexity of the mappings is complicated, since in the
end, the chosen resolution of the image is the dominant factor for the time
required. In
Table \ref{table: countour}, the effects of the polynomial degree and
the chosen tolerance on the overall execution time are described. As a test case, a grid
similar to one in Figure \ref{fig: disk-in-pentagon}, has been computed by using $\sigma = 0.25$
and $\epsilon = 1/100,1/1000,1/10000$, for $p=4,8,12$. Note that for the radial contours the effect of
$\epsilon$ is not as noticeable as for the circular ones due to contours and
 gradients being aligned.

\section{Numerical Experiments}
\begin{figure}[!ht]
\begin{center}
\includegraphics[width=0.15\textwidth]{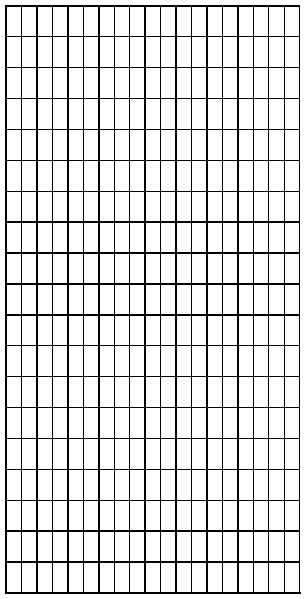} \hspace*{2cm}
\includegraphics[width=0.3\textwidth]{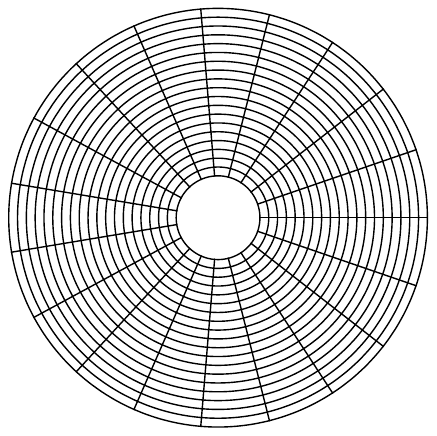}
\caption{Canonical domains $R_h$ and $A_r$ on the left- and right-hand side, respectively. \label{fig: can-domain}}
\end{center}
\end{figure}

\begin{figure}[!ht]
\begin{center}
\includegraphics[width=0.4\textwidth]{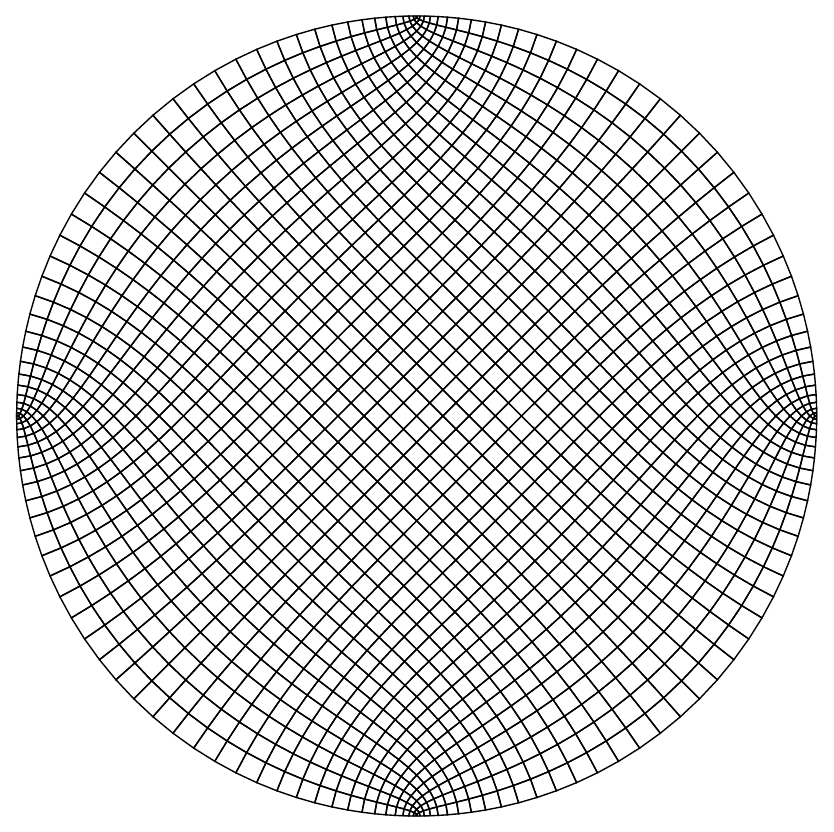} 
\caption{Example of the conformal map of a square onto a disk, first obtained by Schwarz in 1869 \cite{Sch1}. \label{fig: schwarz}  
}
\end{center}
\end{figure}

\begin{figure}[!ht]
\begin{center}
\includegraphics[width=0.4\textwidth]{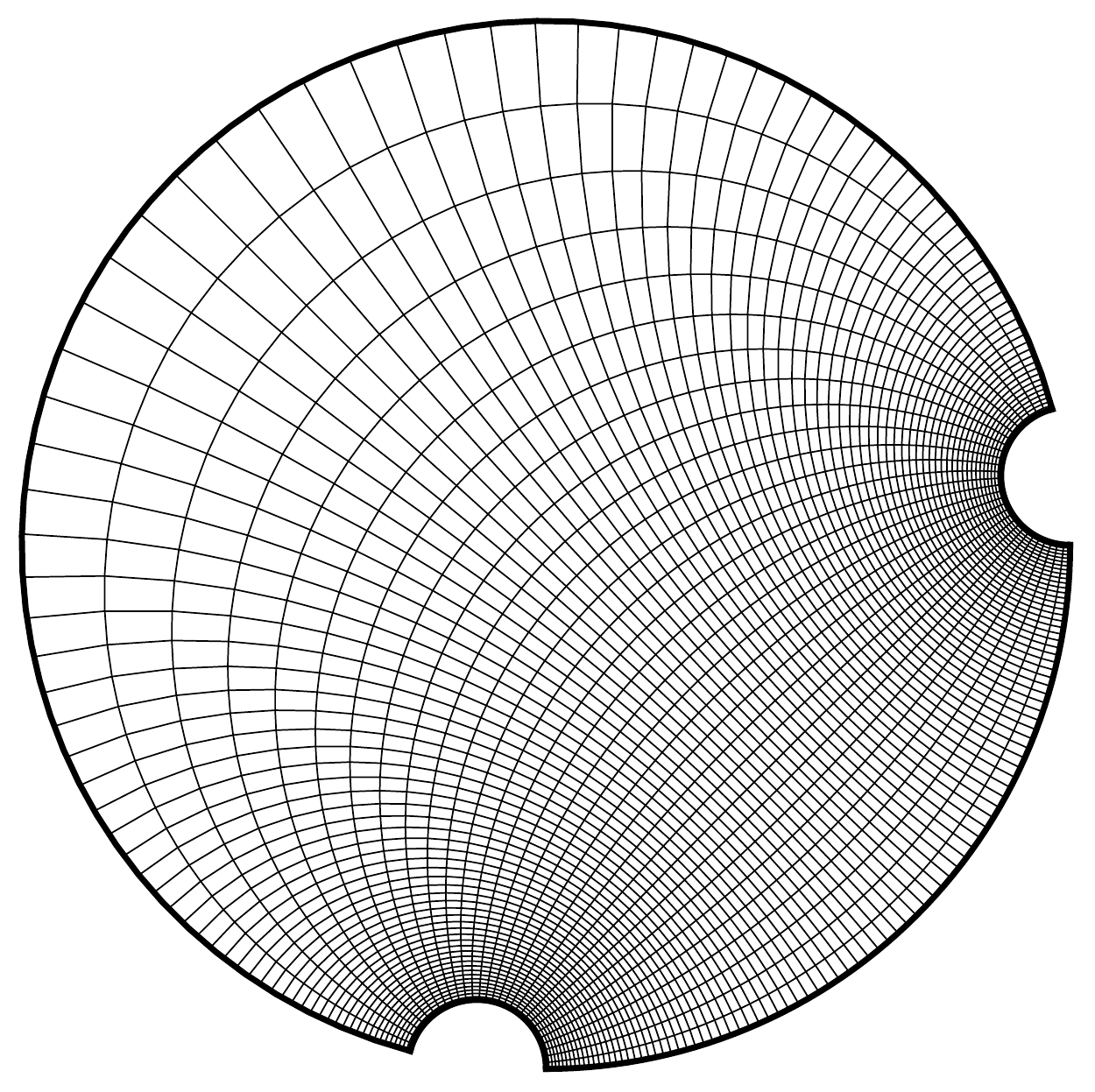}
\caption{The quadrilateral $(\symQuad_A;
  e^{i\pi/12}, e^{i17\pi/12}, e^{i3\pi/2},1)$ and the visualization of
  the pre-image of the rectangular grid  (Figure \ref{fig: can-domain}). \label{fig: quad-a} 
}
\end{center}
\end{figure}

\begin{figure}[!ht]
\begin{center}
\includegraphics[width=0.4\textwidth]{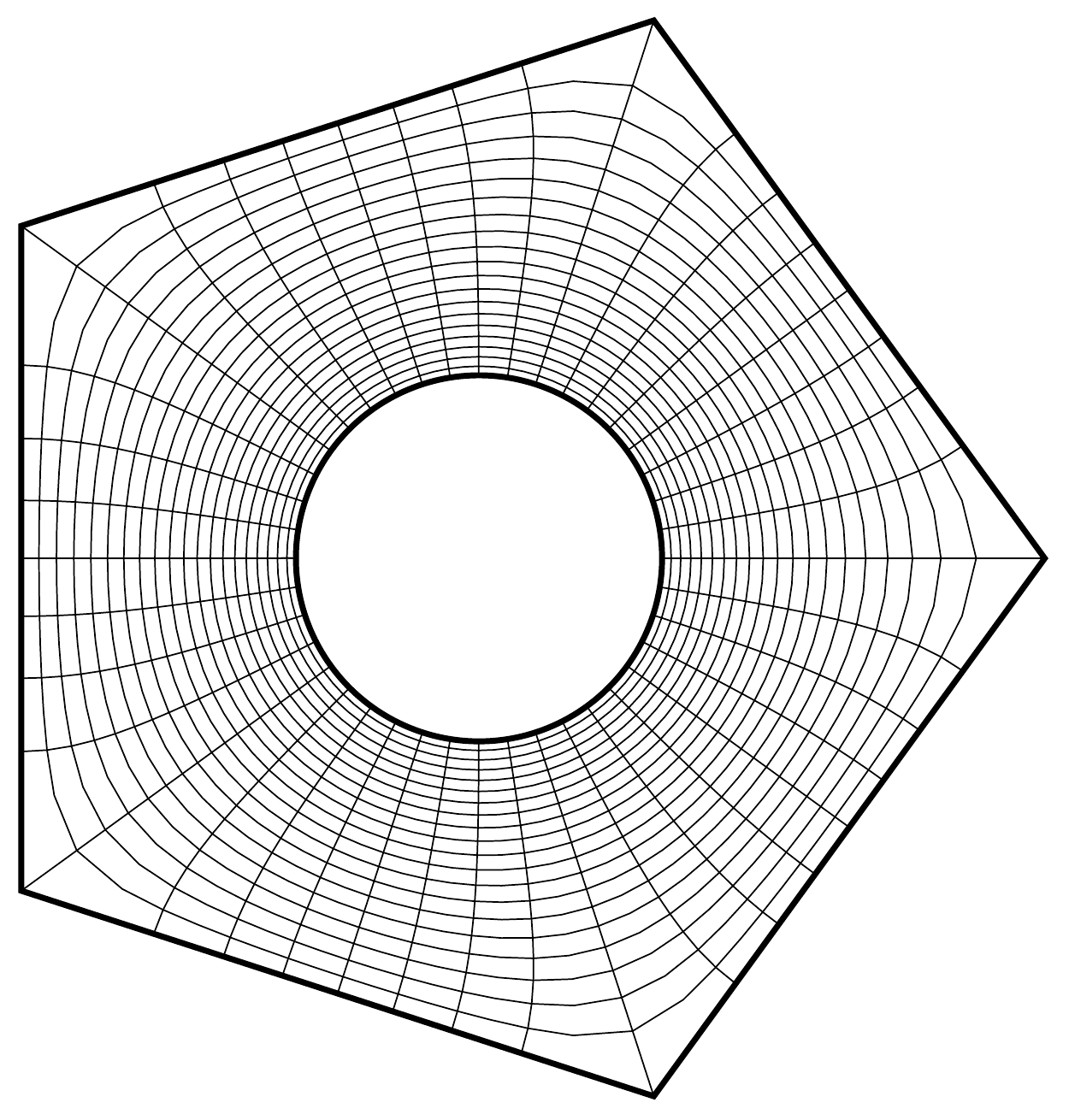}
\caption{Disk in pentagon $(r = 0.4)$ with the pre-image of the
  annular grid  (Figure \ref{fig: can-domain}). \label{fig: disk-in-pentagon} 
}
\end{center}
\end{figure}

Our numerical experiments are divided into three different categories.
First we validate the algorithm against the results obtained by the
Schwarz-Christoffel toolbox and the analytic formula.
Then we study several examples of using our method
to construct conformal mappings from simply (see Figures \ref{fig: schwarz}--\ref{fig: asteroid}) 
or doubly connected (see Figures \ref{fig: cross-in-square}--\ref{fig: droplet}) 
domains onto canonical domains, see Figure \ref{fig: can-domain}, 
with the main
results
summarized in Tables \ref{table: summary-simply} and \ref{table:
  summary-ring}, respectively.

\begin{table}[ht]
\caption{Summary of the tests on simply connected domains. Accuracy is
given as  $\lceil
\log_{10}|1 -\symM(\symQuad) \symM(\symQuadC)|\rceil.$ For the first two 
cases the moduli are known due to symmetry.}\label{table: summary-simply}
\begin{tabular}{|lc|c|c|c|c|}
\hline
Example & ID & $\symM(\symQuad)$ / $\symM(\symQuadC)$ &Accuracy & Figure \\ \hline
Unit Disk & \ref{ex: disk}   &1  / 1 &-13 &
\ref{fig: schwarz}\\
Flower     & \ref{ex: flower} & 1 / 1 & -10 &
\ref{fig: flower}\\
Circular quadrilateral & \ref{ex: circular} & $0.63058735108478$ / &
-13 &   \ref{fig: quad-a}\\
&& $1.585823119159254$ &
&  \\
Asteroid cusp & \ref{ex: asteroidal} & $0.68435408764536$ / &
-9 &  \ref{fig: asteroid}\\
& & $1.4612318657235575$ &
&\\
\hline
\end{tabular}
\end{table}

\begin{table}[ht]
\caption{Summary of the tests on ring domains.  Accuracy is
given as  $\lceil
\log_{10}|1 -\symM(\symQuad) \symM(\symQuadC)|\rceil$, where the
quadrilaterals are the cut domains.}\label{table: summary-ring}
\begin{tabular}{|lc|c|c|c|}
\hline
Example & ID & $\symM(R)$ & Accuracy & Figure  \\ \hline
Disk in regular pentagon & \ref{ex: disk-in-pentagon}  & See Table \ref{table: disk-in-pentagon}. & & \ref{fig: disk-in-pentagon}\\
Cross in square & \ref{ex: cross-in-square} & $0.2862861647287473$ &
-9  & \ref{fig: cross-in-square}\\
Circle in square  & \ref{ex: circle-in-square} & $0.9920378629010557$
& -13 & \ref{fig: disk-in-square}\\
Flower in square & \ref{ex: flower-in-square} & $0.6669554623348065$ &
-8 & \ref{fig: flower-in-square} \\
Circle in L & \ref{ex: circle-in-l} & $1.0935085836560234$ &
-9 & \ref{fig: l-block}\\
Droplet in square & \ref{ex: droplet-in-square} & $0.8979775098918368$
& -9 &  \ref{fig: droplet}\\
\hline
\end{tabular}
\end{table}

\subsection{Setup of the Validation Test} 
Validation of the algorithm for the conformal mapping will be carried out in two cases, first we compare our algorithm with SC Toolbox in a convex and a non-convex quadrilateral. In the second test we parameterized the modulus of a rectangle and map it onto the unit disk.

\begin{figure}[!ht]
\begin{center}
\includegraphics[width=0.37\textwidth]{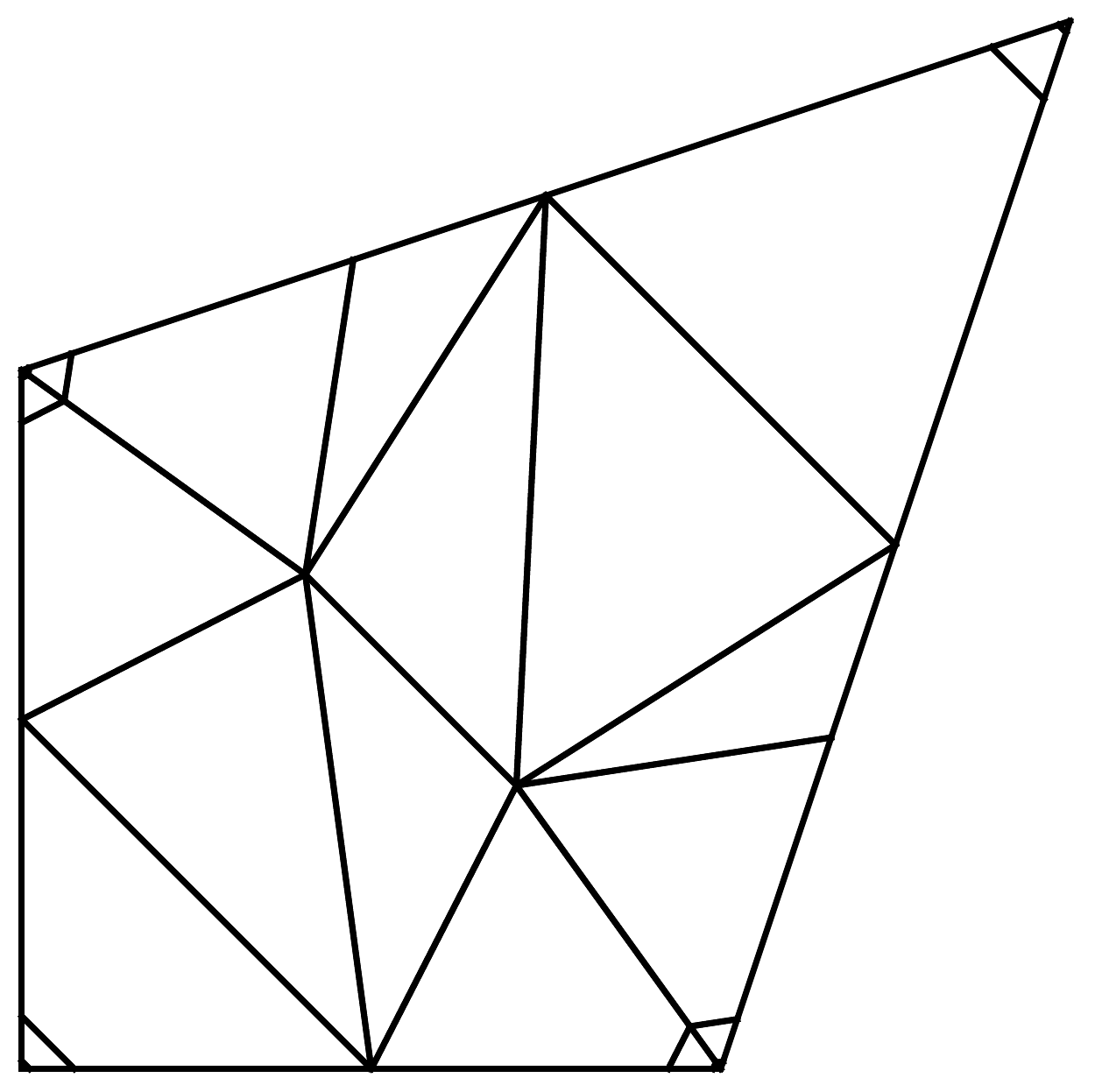} \hspace*{1cm}
\includegraphics[width=0.37\textwidth]{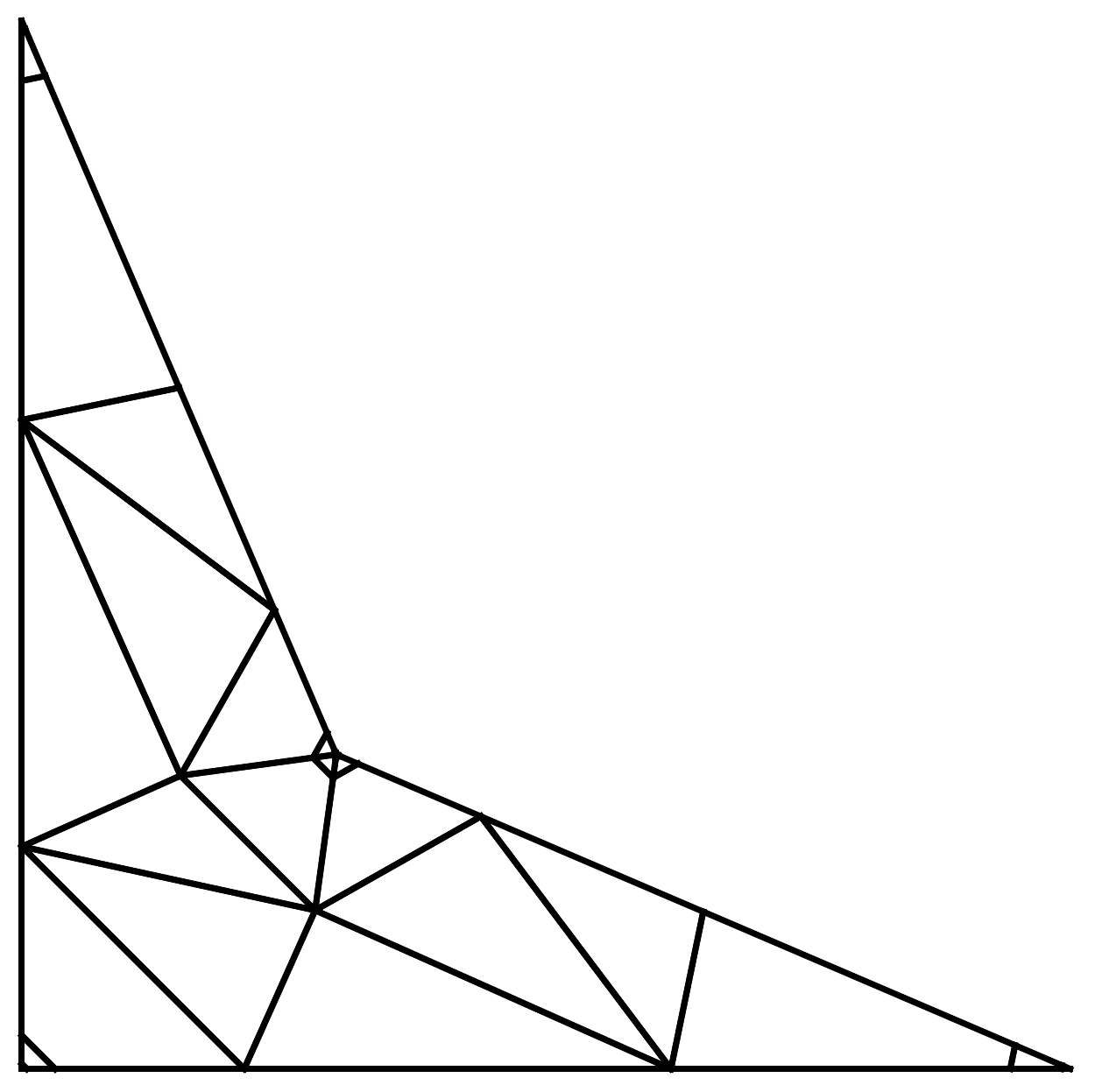}\\
\caption{Geometric mesh of the convex (left-hand side) and the non-convex (right-hand side) quadrilateral used in computing the potential functions. \label{fig: mesh}}
\end{center}
\end{figure}

\begin{figure}[!ht]
\begin{center}
\includegraphics[width=0.45\textwidth]{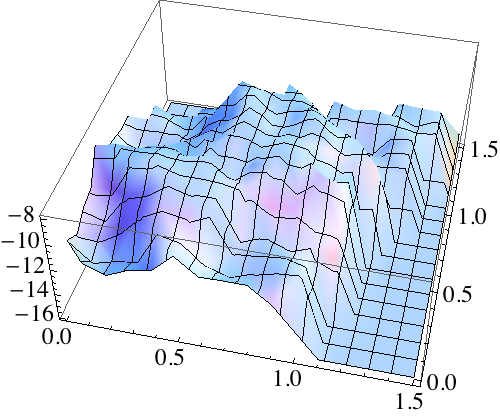} \hspace*{1cm}
\includegraphics[width=0.45\textwidth]{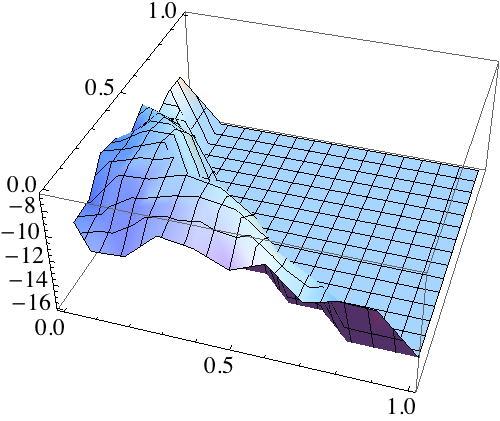}
\caption{Comparison of the convex (left-hand side) and non-convex (right-hand side) quadrilateral between the conjugate function method and SC Toolbox. Results are obtained by taking the logarithm (with base 10) of the test function (\ref{eqn: test-conf}). \label{fig: comparison}}
\end{center}
\end{figure}

The comparison to the SC Toolbox is carried out in the following quadrilaterals: convex quadrilateral $(\symD;0,1,1.5+1.5i,i)$ and non-convex quadrilateral $(\symD;0,1,0.3+0.3i,i)$, and line-segments joining the vertices as the boundary arcs. Then comparisons of the results obtained by the conjugate function method, presented in this paper,
and SC Toolbox by Driscoll \cite{Dri} are carried out. All SC Toolbox tests were carried with the settings \texttt{precision = 1e-14}. Comparison is done by using the following test function
\begin{equation} \label{eqn: test-conf}
\textrm{test}(z) = |f(z) - g(z)|,
\end{equation}
where $f$ and $g$ are obtained by the conjugate function method and SC Toolbox, respectively. The mesh setup of the quadrilaterals and the results of the test function \eqref{eqn: test-conf} are shown in Figure \ref{fig: mesh} and \ref{fig: comparison}, respectively.

All our examples are carried out in the same fashion using the
reciprocal identity \eqref{eqn: recip} and 
a quadrilateral $\symQuad$. 
The test function is
\[
\textrm{rec}(\symQuad) = |\symM(\symQuad) \, \symM(\symQuadC) - 1 |,
\]
which vanishes identically. See also \cite[Section 4]{HRV}.

In the second validation test, we parameterized a rectangle with respect to the modulus $\symM(\symQuad)$ and map the rectangle onto the unit disk. The mapping is given by a composite mapping consisting of a Jacobi's elliptic sine function and a M\"{o}bius transformation.

For every point $(x_j, y_j)$ in the grid on the rectangle $R_h$, where $x_j = j/10$ and $y_j = jh/10$, $j=0,1,2,\ldots,10$, we compute the error $\|e_j\|$ which is simply the Euclidean distance of the image of the initial point $(x_j,y_j)$ computed by the conjugate function method and the analytical mapping.
For a given modulus $\symM(\symQuad)$ the values $\textrm{rec}(\symQuad)$, $\max(\|e_j\|)$, and $\textrm{mean}(\|e_j\|)$,
where the latter two represent the maximal and the mean error over the grid are presented in Table \ref{table: disk-to-rectangle}.

Note that our test function $\textrm{rec}(\symQuad)$ effectively measures the error in energy.
Given the very high accuracy of the results obtained, we are confident that even though
no a priori guarantees for pointwise convergence can be given, the second test is
a valid indication of the global convergence behavior.

\begin{table}[ht]
\caption{The values of $\textrm{rec}(\symQuad)$, $\max(\|e_j\|)$ and  $\textrm{mean}(\|e_j\|)$ for a given $\symM(\symQuad)$.}\label{table: disk-to-rectangle}
\begin{tabular}{|c|c|c|c|}
\hline
$\symM(\symQuad)$ & $\textrm{rec}(\symQuad)$ & $\max(\|e_j\|)$ & $\textrm{mean}(\|e_j\|)$  \\ \hline
$1$   & $8.08242\cdot 10^{-14}$ & $1.87409\cdot 10^{-8}$ & $5.56947\cdot 10^{-10}$\\
$1.2$ & $6.35048\cdot 10^{-14}$ & $7.97889\cdot 10^{-9}$ & $7.49315\cdot 10^{-10}$\\
$1.4$ & $5.52891\cdot 10^{-14}$ & $1.21851\cdot 10^{-8}$ & $6.90329\cdot 10^{-10}$\\
$1.6$ & $8.85958\cdot 10^{-14}$ & $1.10001\cdot 10^{-8}$ & $7.90840\cdot 10^{-10}$\\
$1.8$ & $9.72555\cdot 10^{-14}$ & $1.19005\cdot 10^{-8}$ & $7.31645\cdot 10^{-10}$ \\
$2$   & $9.41469\cdot 10^{-14}$ & $8.56068\cdot 10^{-9}$ & $7.67815\cdot 10^{-10}$ \\
\hline
\end{tabular}
\end{table}

\subsection{Simply Connected Domains} \label{sec: simply-con}
In this section we consider a conformal mapping of a quadrilateral $Q = (\symD; z_1,z_2, z_3, z_4)$ with curved boundaries $\g_1, \g_2, \g_3$, $\g_4$ onto a rectangle $R_h$ such that the vertices $z_1, z_2, z_3, z_4$ maps to $1+ih, ih,0,1$, respectively, and the boundary curves $\g_1, \g_2, \g_3, \g_4$ maps onto the line segments $\g_1', \g_2', \g_3', \g_4'$. We give some examples and applications with illustrations. Simple examples of such domains are domains, where four or more points are connected with circular arcs. Some examples related to numerical methods and the Schwarz-Christoffel formula for such domains can be found in the literature, e.g., \cite{Bro,BP,KP}.


\begin{ex}[Unit disk] \label{ex: disk}
Let $\symD$ be the unit disk. We consider a quadrilateral $\symQuad = (\symD; z_1,z_2, z_3, z_4)$, where $z_j=e^{i\theta_j}$, $\theta_j = (j-1)\pi/ 2$. Note that, because of the symmetry, it follows from \eqref{eqn: recip} that the modulus is $1$. The reciprocal error of the conformal mappings is $4.34 \cdot 10^{-14}$. This example was first given by Schwarz in 1869 \cite{Sch1}.
\end{ex}

\begin{figure}[!ht]
\begin{center}
\includegraphics[width=0.4\textwidth]{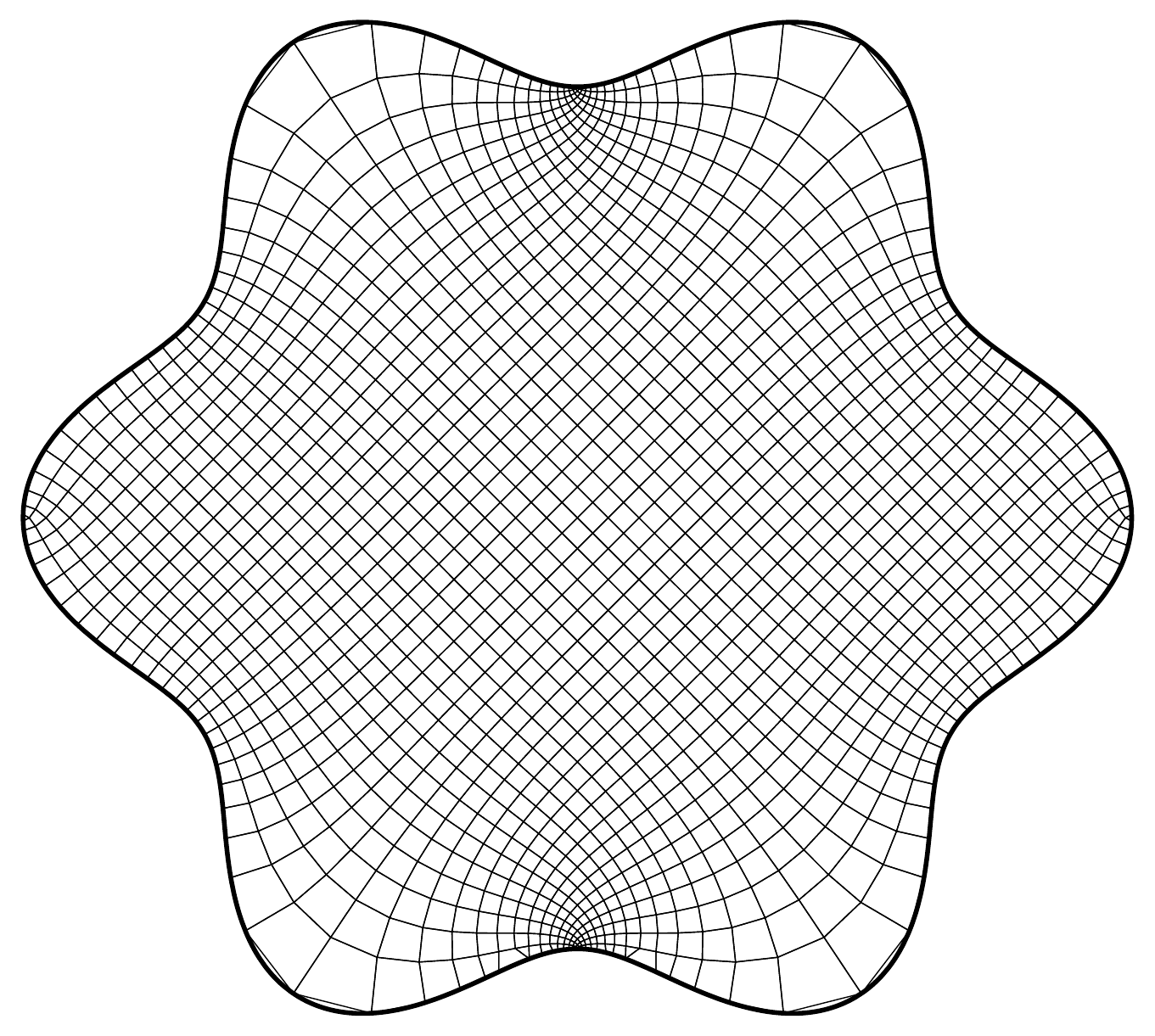}
\caption{Illustration of the flower domain and the
  visualization of the pre-image of the rectangular grid (Figure \ref{fig: can-domain}). \label{fig: flower}  
}
\end{center}
\end{figure}

\begin{ex}[Flower] \label{ex: flower}
Let $\symD$ be the domain bounded by the curve 
\begin{equation} \label{eqn: flower}
r(\theta)= 0.8 + t \cos(n \theta),
\end{equation}
where $0 \le \theta\le 2 \pi$, $n=6$ and $t=0.1$. We consider a quadrilateral $\symQuad = (\symD; z_1,z_2, z_3, z_4)$, where $z_j=r(\theta_j)$, $\theta_j = (j-1)\pi/ 2$; see Figure \ref{fig: flower}. As in Example \ref{ex: disk}, the modulus of $\symQuad$ is 1.
The reciprocal error of the conformal mappings is $3.74 \cdot 10^{-11}$. Several other examples of flower shaped quadrilaterals are given in \cite[Section 8.5]{HRV}.
\end{ex}

\begin{ex}[Circular Quadrilateral] \label{ex: circular}
In \cite{HRV} several experiments of circular quadrilaterals are given. Let us consider a quadrilateral whose sides are circular arcs of intersecting orthogonal circles, i.e.,  angles are $\pi/2$. Let $0<a<b<c<2\pi$ and choose the points $\{1, e^{ia}, e^{ib}, e^{ic}\}$ on the unit circle. Let $\symQuad_A$ stand for the domain which is obtained from the unit disk by cutting away regions bounded by the two orthogonal arcs with endpoints $\{ 1, e^{ia} \}$ and $\{ e^{ib}, e^{ic} \}$, respectively. Then $\symQuad_A$ determines a quadrilateral $(Q_A ; e^{ia}, e^{ib}, e^{ic},1)$. Then for the triple $(a,b,c) = (\pi/12, 17\pi/12,3\pi/2)$, the modulus $\symM(\symQuad_A) =0.630587351084775$ and $\symM(\symQuadC_A) = 1.5858231191592544$. The reciprocal error of the conformal mapping is $1.68 \cdot 10^{-13}$. 
\end{ex}

\begin{figure}[!ht]
\begin{center}
\includegraphics[width=0.4\textwidth]{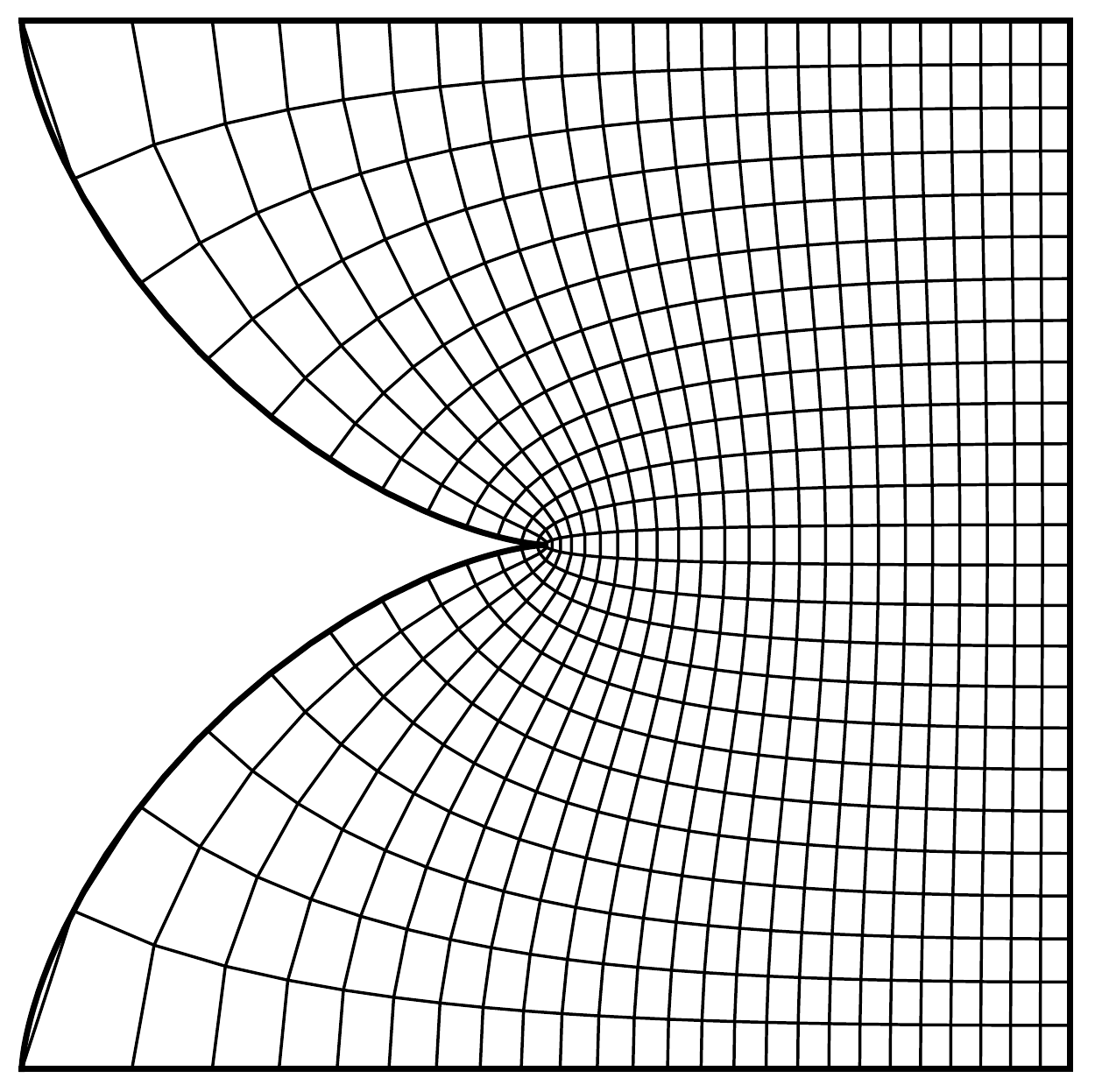}
\caption{Asteroid cusp domain with the pre-image
  of the rectangular grid  (Figure \ref{fig: can-domain}). \label{fig: asteroid}  
}
\end{center}
\end{figure}

\begin{ex}[Asteroid Cusp] \label{ex: asteroidal}
Asteroid cusp is a domain $\symD$ given by a
\begin{equation} \label{eqn: box}
G_c=\{(x,y):|x|<c,\,|y|<c\},
\end{equation} 
where $c=1$ and the left-hand side vertical boundary line-segment is replaced by the following curve
\[
r(t) = -1 + \cos^3 t + i\sin^3 t,\ t \in [-\pi/2,\pi/2].
\]
We consider a quadrilateral $\symQuad = (\symD; 1-i,\, 1+i, \, -1+i,
\, -1-i)$. The reciprocal error of the conformal mappings is of the
order $10^{-10}$. The modulus $\symM(\symQuad) = 0.68435408764536$ and $\symM(\symQuadC) = 1.4612318657235575$. The domain is illustrated in Figure \ref{fig: asteroid}.
\end{ex}

\subsection{Ring Domains}

In this section we shall give several examples of conformal mapping from a ring domain $R$ onto an annulus $A_r$. It is also possible to use the Schwarz-Christoffel method, see \cite{Hu}. For symmetrical ring domains a conformal mapping can be obtained by using Schwarz' symmetries.

\begin{ex}[Disk in Regular Pentagon] \label{ex: disk-in-pentagon}
Let $\symD$ be a regular pentagon centered at the origin and having short radius (apothem) 
equal to $1$ such that the corners of the pentagon are $z_k = \sec(\pi/5) e^{2\pi ik/5}$, $k=0,1,2,3,4$. Let $\symDisk(r)=\{z \in \C:|z|\leq r\}$. We consider a ring domain 
$R = \symD \backslash \symDisk(r)$ and compute the modulus $\symM(R)$ and the exponential of the modulus $e^{\symM(R)}$. Results are reported in Table \ref{table: disk-in-pentagon} with the values $e^{\symM(R)}$ from \cite[Example 5]{BSV} in the fourth column.

\begin{center}
\begin{table}[ht]
\caption{The values $\symM(R)$ and $e^{\symM(R)}$.}\label{table: disk-in-pentagon}
\begin{tabular}{|c|c|c|c|}
\hline
$r$ & $\symM(R)$ & $\exp(\symM(R))$ & \cite[Example 5]{BSV} \\
\hline
$0.1$ & $2.35372035858745$ & $10.524652459913115$ & $10.5246525$ \\
$0.4$ & $0.9674246001764809$ & $2.631159438480101$ & $2.631159439$ \\
$0.9$ & $0.15070188000332954$ & $1.1626499971978235$ & $1.1626499972$ \\
$0.99$ & $0.03276861064365647$ & $1.0333114143138304$ & $1.03331141431$ \\
$0.999$ & $0.00934656029871744$ & $1.0093903757950962$ & $1.00939037579$ \\
\hline
\end{tabular}
\end{table}
\end{center}
\end{ex}

\begin{figure}[!ht]
\begin{center}
\includegraphics[width=0.4\textwidth]{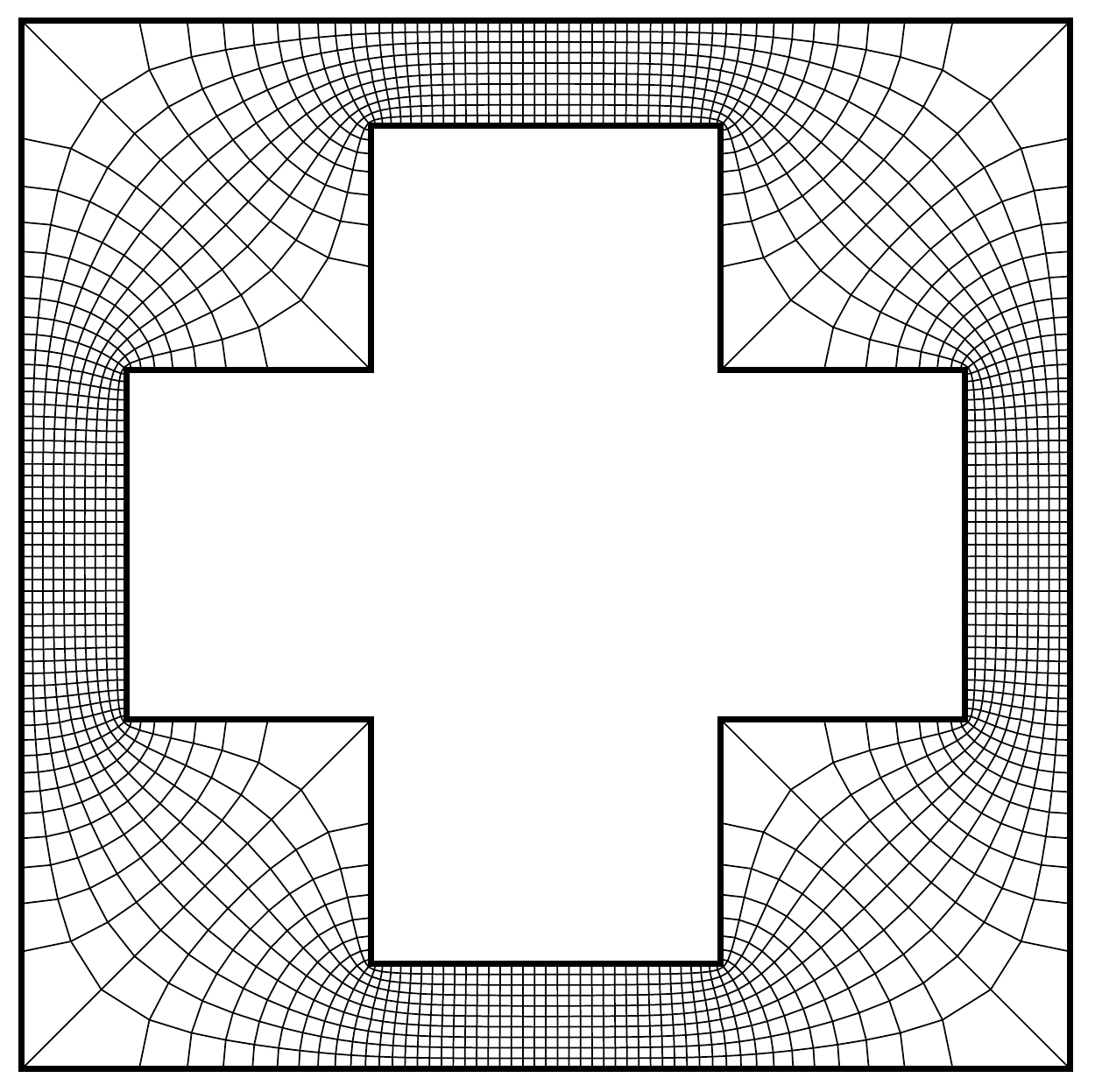}
\caption{The ring domain $G_c \backslash
  G_{ab}$, where $a=0.5, b=1.2, c=1.5$, with the pre-image of the
  annular grid  (Figure \ref{fig: can-domain}).  \label{fig: cross-in-square} 
}
\end{center}
\end{figure}

\begin{ex}[Cross in Square] \label{ex: cross-in-square}
Let $G_{ab}=\{(x,y):|x|\leq a,\, |y|\leq b\}\cup \{(x,y):|x|\leq b,\, |y|\leq a\}$, and $G_c$ as in (\ref{eqn: box}), where $a<c$ and 
$b<c$. Then the domain cross in square is a ring domain $R = G_c \backslash G_{ab}$, see Figure \ref{fig: cross-in-square}. 
The reciprocal error of the conformal mapping is of the order $10^{-10}$. The modulus $\symM(R) = 0.2862861647287473$.
\end{ex}

\begin{figure}[!ht]
\begin{center}
\includegraphics[width=0.4\textwidth]{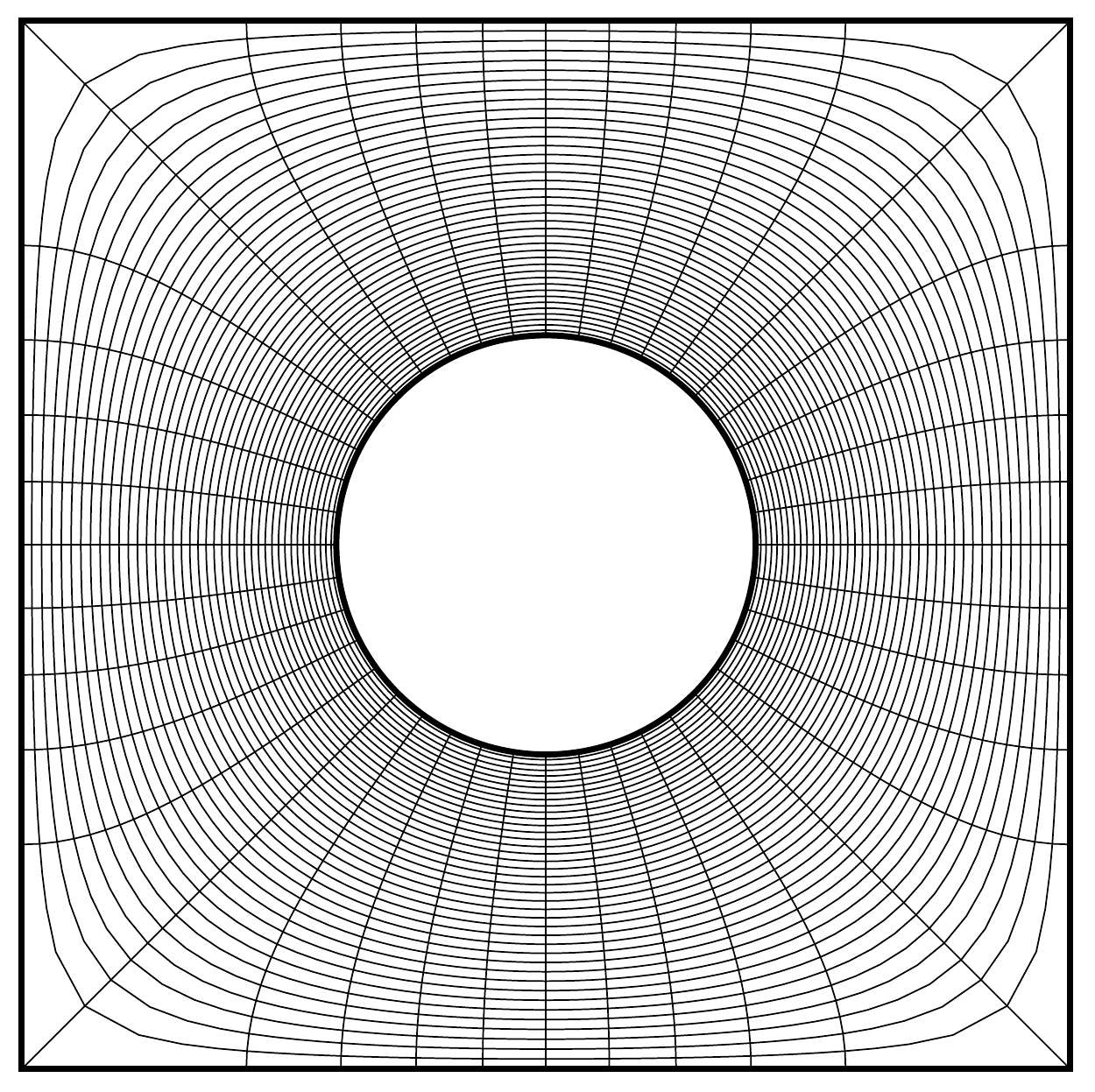}
\caption{Disk in a square domain  with the pre-image of the
  annular grid  (Figure \ref{fig: can-domain}). \label{fig: disk-in-square} 
}
\end{center}
\end{figure}

\begin{ex}[Circle in Square] \label{ex: circle-in-square}
Let $\symD$ be the unit disk. Then we consider a ring domain $R = G_c \backslash \symD$, 
where $c=1.5$, see Figure \ref{fig: disk-in-square}. The reciprocal error of the conformal 
mapping is of the order $10^{-14}$. The modulus $\symM(R) = 0.9920378629010557$.
\end{ex}

\begin{figure}[!ht]
\begin{center}
\includegraphics[width=0.4\textwidth]{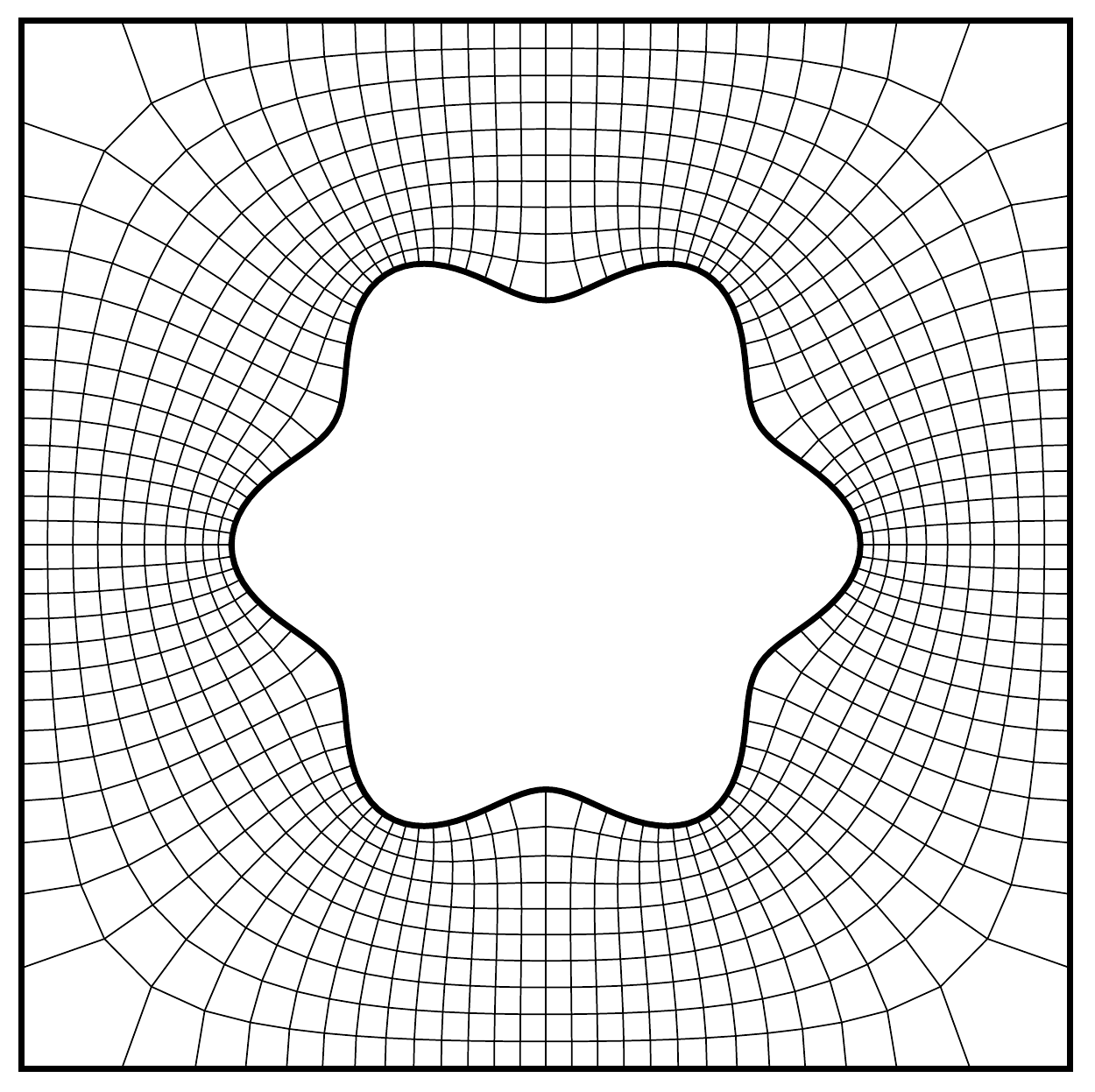}
\caption{Flower in a square domain with the pre-image of the
  annular grid  (Figure \ref{fig: can-domain}). \label{fig: flower-in-square} 
}
\end{center}
\end{figure}

\begin{ex}[Flower in Square] \label{ex: flower-in-square}
Let $\symD$ be a domain bounded by the curve (\ref{eqn: flower}). 
Then we consider a ring domain $R = G_c \backslash \symD$, where $G_c$ is 
given by (\ref{eqn: box}) and $c=1.5$. See Figure \ref{fig: flower-in-square} for 
the illustration. The reciprocal error of the conformal mapping is of the order $10^{-9}$. 
The modulus $\symM(R) = 0.6669554623348065$.
\end{ex}

\begin{figure}[!ht]
\begin{center}
\includegraphics[width=0.5\textwidth]{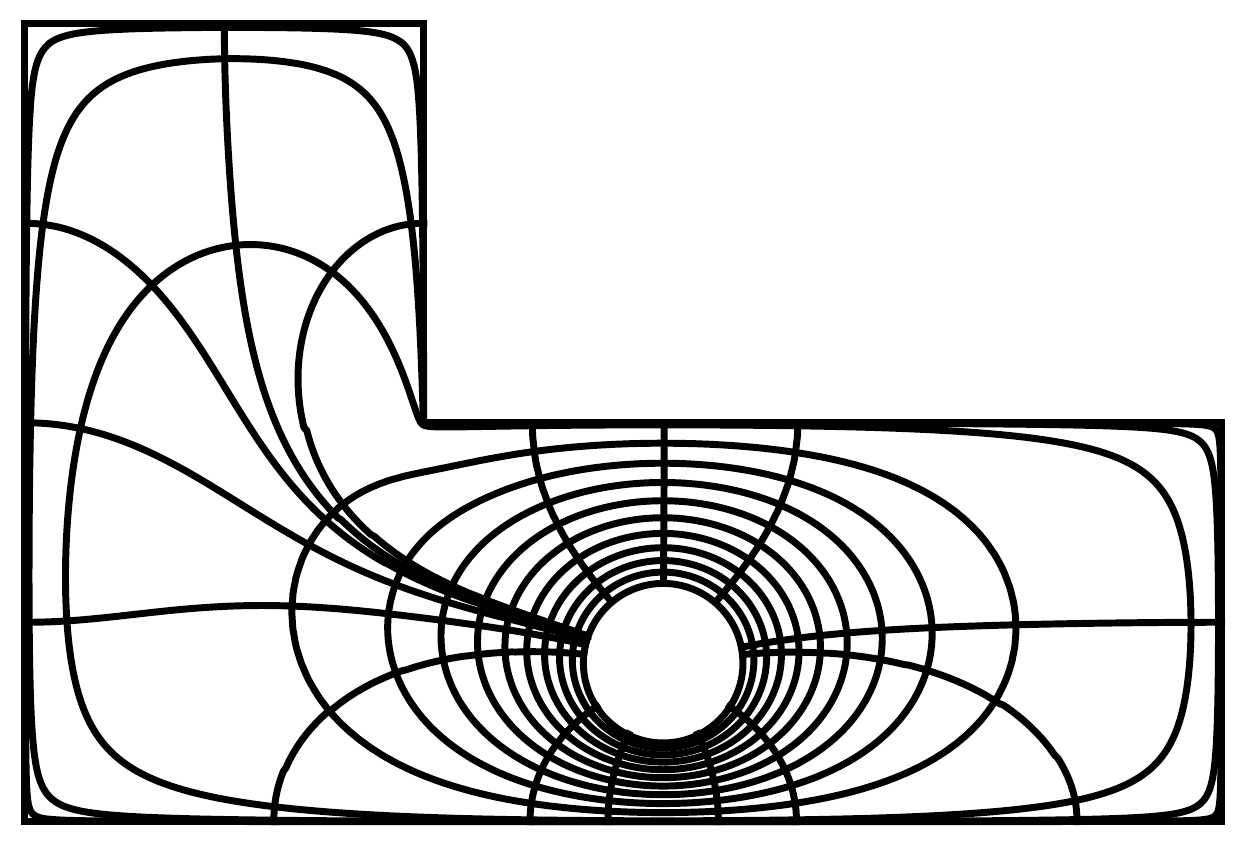} \vspace*{0.25cm}

\caption{L-shaped domain with a circular hole with the pre-image of the
 non-uniform  annular grid of Example \ref{ex: circle-in-l}. \label{fig: l-block} 
}
\end{center}
\end{figure}

\begin{ex}[Circle in L] \label{ex: circle-in-l}
Let $L_1=\{ z\in \C : 0 < \textrm{Re}(z) < a,\, 0<\textrm{Im}(z) <b \}$ and 
$L_2= \{ z\in {\mathbb C} : 0 < \textrm{Re}(z) < d,\, 0<\textrm{Im}(z) <c \}$,
where $0<d<a, \; 0<b<c$. Then $L(a,b,c,d) =L_1\cup L_2$ is called an L-domain. 
Suppose that $\symDisk(z_0,r) = \{ z \in \C : |z-z_0| < r \}$.  
We consider a ring domain $R = L(a,b,c,d) \backslash \symDisk(z_0,r)$, 
where $(a,b,c,d) = (3, 1, 2, 1)$, $z_0 = 8/5 + 2i/5$, and $r = 1/5$. 
See Figure \ref{fig: l-block}.
%

In order to better illustrate the details of the mapping, a non-uniform 
grid has been used. For the real component the points $x$ are
\[
x = \{k / 10\ : \ k=0,1,\ldots,9\} \cup \{99/100,999/10000,9999/10000,1\}.
\]
For the imaginary component the points $y$ are chosen on purely aesthetic basis as:
\begin{align*}
y = \{k / 10\ : \ k=1,2,\ldots,9\} &\ \cup\\ \{0.316225, 0.324008, & 0.327831, 0.329278, 0.331005, 0.687482\}.
\end{align*}

The reciprocal error of the conformal mapping is of the order $10^{-10}$. 
The modulus $\symM(R) = 1.0935085836560234$.
\end{ex}

\begin{figure}[!ht]
\begin{center}
\includegraphics[width=0.4\textwidth]{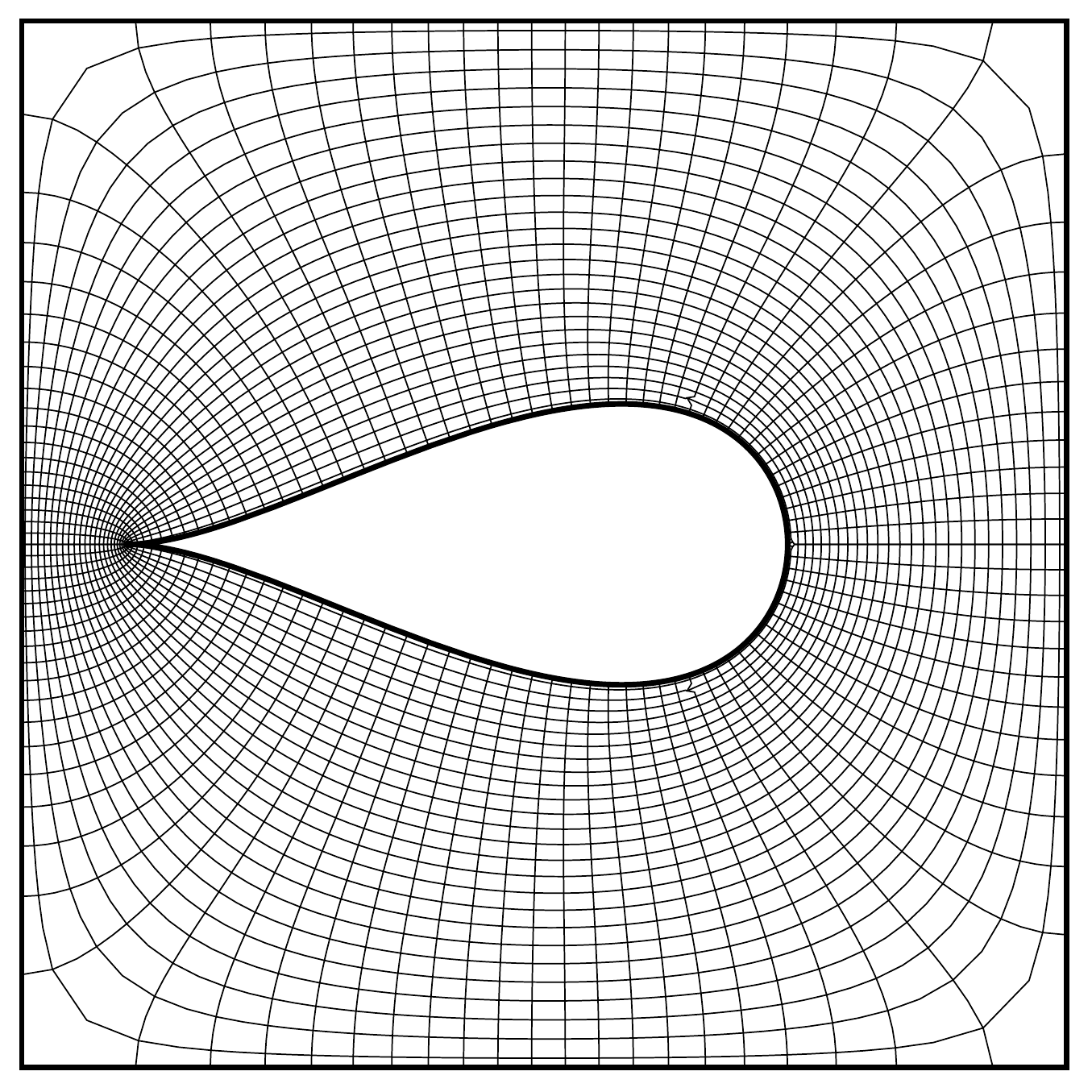}
\caption{Droplet in square with the pre-image of the
  annular grid  (Figure \ref{fig: can-domain}). \label{fig: droplet} 
}
\end{center}
\end{figure}

\begin{ex}[Droplet in Square] \label{ex: droplet-in-square}
Let $Q_D$ be bounded by a Bezier curve:
\[
r(t) = \frac{1}{640} \left(45 t^6+75 t^4-525 t^2+469\right) +
\frac{15}{32} t
   \left(t^2-1\right)^2 i,\ t \in [-1,1].
\]
Then the domain droplet in square is a ring domain $R = G_c \backslash Q_D$, where $G_c$ in given in the first example concerning ring domains. For visualization, see Figure \ref{fig: droplet}. The reciprocal error of the conformal mapping is of the order $10^{-10}$. 
The modulus $\symM(R) = 0.8979775098918368$.
\end{ex}

\begin{acknowledgement}
We thank T.A.~Driscoll, R.M.~Porter and M.~Vuorinen for their valuable comments on this paper.\end{acknowledgement}

\end{document}